\documentclass[11pt]{article}

\usepackage{amssymb,epsfig}
\usepackage[mathscr]{euscript}

\def\lf{\hfill \break}
\def\ds{\displaystyle}

\newtheorem{proposition}{Proposition}
\newtheorem{remark}{Remark}
\newtheorem{definition}{Definition}
\newtheorem{lemma}{Lemma}
\newtheorem{theorem}{Theorem}
\newtheorem{corollary}{Corollary}

\newcommand{\goth}[1]{\ensuremath{\mathfrak{#1}}}
\newcommand{\bbox}{\normalsize {}%
        \nolinebreak \hfill $\blacksquare$ \medbreak \par}

\newcommand{\R}[1][]{\ensuremath{{\mathbb{R}^{#1}} }}
\newcommand{\C}[1][]{\ensuremath{{\mathbb{C}^{#1}} }}
\newcommand{\CP}[1][]{\ensuremath{{\mathbb{CP}^{#1}} }}
\newcommand{\N}[1][]{\ensuremath{{\mathbb{N}^{#1}} }}
\newcommand{\Z}[1][]{\ensuremath{{\mathbb{Z}^{#1}} }}

\newcommand{\T}[1][]{\ensuremath{{\mathbb{T}^{#1}} }}

\newcommand{\G}{\ensuremath{\mathcal{G}}}
\newcommand{\GC}{\ensuremath{\mathcal{G}^{\C}}}
\newcommand{\B}{\ensuremath{\mathcal{B}}}
\newcommand{\g}{\ensuremath{\goth{g}}}
\newcommand{\gC}{\ensuremath{\goth{g}^{\C}}}
\renewcommand{\b}{\ensuremath{\goth{b}}}

\newcommand{\beps}{\bar{\epsilon}}
\newcommand{\bz}{\bar{z}}

\newcommand{\ad}{\mathrm{ad}}

\def\<{\langle} \def\>{\rangle}
\newcommand{\pa}[3][]{ \frac{ \partial^{#1} {#2} }{ \partial {#3}^{#1} }  }
\newcommand{\1}{\mathrm{1 \hspace{-0.25em} l}}
\renewcommand{\Re}{\mathrm{Re}}
\renewcommand{\Im}{\mathrm{Im}}
\newcommand{\trsp}[1]{{}^{t}\!#1}
\newcommand{\Gg}{\ensuremath{\goth{G}}}

\title{Hamiltonian stationary Lagrangian surfaces in \C[2]}
\author{Fr\'ed\'eric H\'elein, Pascal Romon}
\date{28-5-1999}

\begin{document}
\maketitle

\begin{abstract}
We study Hamiltonian stationary Lagrangian surfaces in \C[2], i.e.
Lagrangian surfaces in \C[2] which are stationary points of the
area functional under smooth Hamiltonian variations.
Using loop groups, we propose a formulation of the equation as a completely 
integrable system. We construct a Weierstrass type representation and
produce all tori through either the integrable systems machinery or 
more direct arguments. \lf
\lf
\end{abstract}

\lf

\section{Introduction}

This paper addresses the study of Hamiltonian stationary oriented
Lagrangian surfaces in a symplectic Euclidean vector space of dimension
4, using techniques of completely integrable systems. The ambient space may
be seen as $\C^2$ with, using complex coordinates $z^1=x^1+iy^1$ and
$z^2=x^2+iy^2$, the symplectic form
$\omega = dx^1\wedge dy^1 + dx^2\wedge dy^2$ and the canonical scalar
product. The Lagrangian surfaces in $\C^2$ are the immersed surfaces 
on which  the restriction of $\omega$ vanishes. On the set of oriented
Lagrangian surfaces $\Sigma$ in $\C^2$, we let the area functional to 
be

$${\cal A}(\Sigma ) = \int _{\Sigma}dv,$$
where the volume form $dv$ is defined using the induced metric on $\Sigma$. A
critical point of this functional is a Lagrangian surface such that $\delta
{\cal A}(\Sigma )(X)=0$ for any compactly supported smooth vector field $X$
on $\C^2$, satisfying some particular constraint: if $X$ is arbitrary 
we just say that $\Sigma$ is {\em stationary} (it is actually a minimal
surface in $\C^2\simeq \R[4]$), if $X$ is Lagrangian, i. e. its 
flow preserves
Lagrangian surfaces, $\Sigma$ is called {\em Lagrangian 
stationary}\footnote{called isotropic minimal in \cite{CM}}, and
lastly if $X$ is Hamiltonian, i. e.
$X = -J\nabla h = {\partial h\over \partial y^1} {\partial \over \partial 
x^1}
- {\partial h\over \partial x^1} {\partial \over \partial y^1}
+ {\partial h\over \partial y^2} {\partial \over \partial x^2}
- {\partial h\over \partial x^2} {\partial \over \partial y^2}$, for
some $h\in {\cal C}^{\infty}_c(\C^2, \R)$, $\Sigma$ is called 
{\em
Hamiltonian stationary}\footnote{called E-minimal in \cite{CM} 
and H-minimal in \cite{O2} and \cite{CU}}.

The first variation of the area involves the Lagrangian angle: if $m$ is a
point in $\Sigma$ and if $(e_1,e_3)$ is a direct orthonormal basis of
$T_m\Sigma$, $dz^1\wedge dz^2(e_1,e_3)$ is a complex number of modulus
equal to 1, which we can denote $e^{i\beta}$, for some real number
$\beta$. It builds up a map $\beta$ from $\Sigma$ to $\mathbb{R}/2\pi 
\mathbb{Z}$.
This map is a part of the full Gauss map of the immersion of $\Sigma$. The
mean curvature vector $H$ on $\Sigma$ is then given by $H = J\nabla \beta$,
and thus

$$\delta {\cal A}(\Sigma )(X) = \int _{\Sigma} \langle X, H\rangle dv =
 \int _{\Sigma} \langle -J\nabla h, J\nabla \beta \rangle dv =
\int _{\Sigma} \langle \nabla h, \nabla \beta \rangle dv,$$
see \cite{O1} for more details. Hence Hamiltonian stationary surfaces are
characterized by the equation $-\Delta \beta = 0$,
where $\Delta$ is the
Laplace operator on $\Sigma$, which comes from the induced metric. Surfaces
such that $\beta$ is constant (or $H=0$) are a particular case, called {\em
special Lagrangian surfaces} by R. Harvey and H.B. Lawson \cite{HaL}: they are 
actually
area minimizing since calibrated by $e^{-i\beta}dz^1\wedge dz^2$.

Examples of Hamiltonian stationary surfaces are the standard 
square tori
$T_r = \{(z^1,z^2)\in \C^2/ |z^1| = |z^2| = r\sqrt{2}\}$ and
``rectangular'' variants
$T_{r_1,r_2} = \{(z^1,z^2)\in \C^2/ |z^1|/r_1 = |z^2|/r_2 = 
\sqrt{2}\}$.
These are candidates to be area minimizing with respect to Hamiltonian
deformations as conjectured by Y.G. Oh \cite{O1, O2}. More recently, 
I. Castro and F. Urbano \cite{CU} have constructed more exotic examples 
of Hamiltonian stationary
tori. Beside these explicit instances, R. Schoen and J. Wolfson announce
various existence and partial regularity results and in particular a proof 
of
the existence of a smooth solution to the Plateau problem in
$\C^2$ \cite{ScW}.

A motivation to study Hamiltonian stationary surfaces is for instance the
following model of incompressible elasticity. If $(\phi ,\psi )$ is a
diffeomorphism between two two-dimensional domains $U$ and $U'$, which is
incompressible, i. e. 
${\partial \phi \over \partial x}{\partial \psi \over \partial y}
- {\partial \phi \over \partial y}{\partial \psi \over \partial x} = 1$
everywhere,
and which minimizes the area of the graph functional
$\int _U\sqrt{2 + |\nabla \phi |^2 + |\nabla \psi |^2}dxdy$ among
all possible incompressible diffeomorphisms with the same boundary data,
then its
graph
$\Sigma = \{ (x,y,\phi (x,y),-\psi (x,y))/(x,y)\in U\}$ is Hamiltonian
stationary Lagrangian and conversely. Such a problem has been considered 
by J. Wolfson in \cite{W}.
Also Hamiltonian stationary surfaces offer a nice generalization of the
minimal surface theory. The conjecture of Y.G. Oh above is an interesting
generalisation of isoperimetric inequalities.
Such an inequality would be related to many questions in symplectic 
geometry, as
illustrated by C. Viterbo \cite{V}, who also gave a lower bound for the area
functional of a torus. Also special Lagrangian surfaces has appeared in
recent developments in mathematical Physics, in M-theory \cite{AFS}, and about
Mirror symmetry for Calabi-Yau manifolds: see for example \cite{SYZ}, where A. 
Strominger, S.T. Yau
and E. Zaslow proposed that the moduli space of special Lagrangian
surfaces in a Calabi-Yau is related to the mirror of the manifold.

Our aim here is to show that the set of Hamiltonian stationary Lagrangian 
surfaces in
$\C^2$
forms a completely integrable system, and to use ideas from the
Adler-Kostant-Symes theory in a
similar way as it was done by F. Burstall, D. Ferus, F. Pedit, U. Pinkall
\cite{BFPP} and
J. Dorfmeister, F. Pedit, H. Wu \cite{DPW} for harmonic maps between a surface
and a homogeneous manifold, or by F. H\'elein \cite{H2} for Willmore surfaces.
(See also \cite{U, SWi, Hi, FP, FPPS, DH} about previous results.)
Our main results are: a formulation of the Hamiltonian stationary surfaces
problem in terms of a family depending on a complex parameter of curvature
free connections (a characteristic feature in integrable systems); a
correspondance between conformal immersions of Hamiltonian stationary
surfaces in $\C^2$ and holomorphic maps into $\C^3$ (similar to \cite{DPW}); 
a proof that all Hamiltonian stationary tori in $\C^2$ are obtained by 
a finite type construction (this is similar to \cite{BFPP}); lastly 
a construction of all such tori by integrating linear elliptic equations.

From the point of view of the theory of completely integrable systems, we
obtain an original (at least for us!) example of situation where:
\begin{itemize}
\item the family of curvature free connections has the form
$\alpha _{\lambda} = \lambda ^{-2}\alpha '_2 + \lambda ^{-1}\alpha '_{-1} +
\alpha _0 + \lambda \alpha ''_1 + \lambda ^2\alpha ''_2$
instead of $\lambda ^{-1}\alpha '_1 + \alpha _0 + \lambda \alpha ''_1$ as
in many integrable systems,
\item the situation is almost linear and, in some situations, simplifies in
such a way that we could present the results without these techniques.
\end{itemize}
However we choose to expose the full machinery in our situation since this
is the way we obtained all the constructions here and it seems to 
illuminate
how completely integrable systems work.\\

Our paper is organized as follows.
In section 2 we present the symmetry
group of affine isometries of $\R[4]$ preserving the symplectic form
and the modelisation of conformal immersions of Hamiltonian stationary
Lagrangian surfaces using moving frames. A Cartan decomposition of the Lie
algebra appears to be the key of the formulation.
In section 3
we show that the construction of  conformal immersions of
Hamiltonian stationary simply connected surfaces is equivalent to solving 
three simple linear
PDE's as follows:
let $\beta$ to be a real harmonic map on a simply connected domain 
$\Omega$; we solve on $\Omega$ the
linear equation

$${\partial u\over \partial \bar{z}} = {1\over 2}{\partial \beta 
\over \partial z}J\bar{u},$$
for $u=\trsp{(a/2, b/2, -ia/2,ib/2)}$ and $a$ and $b$ complex valued
functions. Then we integrate the equation
\[
dX = e^{\beta J/2}(udz + \bar{u}d\bar{z})
\]
to obtain a map $X$ to $\R[4]$. Then $X$ is a weakly conformal
Hamiltonian stationary Lagrangian map.
We use these ideas to deduce explicit parametrizations of all tori and we 
identify known examples: the standard torus and the surfaces of I. Castro
and F. Urbano and we show other examples. 
In section 4 we introduce
loop groups and twisted loop groups and we prove various Riemann-Hilbert 
and
Birkhoff-Grothendieck decomposition results.
In section 5 we use the
previous results to establish a Weierstrass type representation.
In section 6 we use the finite gap ideas in integrable systems and prove 
that, for Hamiltonian stationary conformal
immersions of tori, the set of solutions splits
into a countable union of vector spaces (``finite type'' solutions). 
Lastly we link this formulation with the
one in section 3.

We point out that our results could be generalized to Hamiltonian
stationary Lagrangian conformal immersions in $\CP^2$ (or isotropic 
surfaces in higher dimensional K\"ahler
manifolds). This will be the subject of a forthcoming paper.

{\em Aknowledgements}: the authors wish to thank Franz Pedit for his 
remarks during the preparation of this paper.

\section{Moving frames and groups}

\subsection{Symmetry groups for symplectic Euclidean affine 4-spaces}

Let $E^4$ be an affine oriented Euclidean symplectic space and $\vec{E}^4$ 
the
associated oriented Euclidean vector space. We denote by $\langle 
.,.\rangle$ the scalar
product and $\omega$ the symplectic form on $\vec{E}^4$. There exists a 
unique complex structure $J$ on $\vec{E}^4$, such that
$\omega (x,y) = \langle Jx, y\rangle$, $\forall x,y \in \vec{E}^4$.
We denote by ${\cal F}$, the set of all orthonormal 
bases $e = (e_1,e_2,e_3,e_4)$ of $\vec{E}^4$, such that $e_2 = Je_1$ and 
$e_4 = Je_3$. We choose an origin $O$ in $E^4$ and an orthonormal basis of $\vec{E}^4$, 
$(\epsilon _1, \epsilon _2, \epsilon _3, \epsilon _4)\in {\cal F}$. In
the corresponding coordinates $( x^1, x^2, x^3, x^4)$, the symplectic form
reads
\[
\omega = dx^1\wedge dx^2 + dx^3\wedge dx^4.
\]
And the complex structure $J$ has the matrix
\[
L_i = \left( \begin{array}{cccc}
0 & -1 & 0 & 0\\
1 & 0 & 0 & 0\\
0 & 0 & 0 & -1\\
0 & 0 & 1 & 0
\end{array} \right)
\]
(the meaning of that notation will become clear below).

The relevant symmetry groups here are
\begin{itemize}
\item ${\cal G}$, the group of affine transformations of $E^4$ which
preserve $\langle .,.\rangle$ and $\omega$ (or alternatively which 
preserve $\langle .,.\rangle$ and $J$)
\item $\vec{\cal G}$, the group of linear transformations of
$\vec{E}^4$ which preserve $\langle .,.\rangle$ and $\omega$ (or $\langle 
.,.\rangle$ and $J$), which we
may view as a subgroup of ${\cal G}$, namely the isotropy group at 0.
\end{itemize}

Let us analyze first $\vec{\cal G}$. A first description of $\vec{\cal G}$ 
is obtained by the
identification of $\vec{E}^4$ through the quaternions $\mathbb{H}$:

$$\begin{array}{cccl}
Q: & \vec{E}^4\simeq \R[4] & \longrightarrow & \mathbb{H}\\
 & x^1\epsilon _1+x^2\epsilon _2+x^3\epsilon _3+x^4\epsilon _4\simeq (x^1, 
x^2, x^3, x^4) & \longmapsto &  x^1+ ix^2+jx^3+kx^4.
\end{array}$$
Let $S^3_{\mathbb{H}}=\{p\in \mathbb{H}/|p|=1\}$. To each pair $(p,q)\in
S^3_{\mathbb{H}}\times S^3_{\mathbb{H}}$ corresponds a rotation $G_{(p,q)}\in
SO(4)$ defined by: $\forall x\in \R[4]$,

$$Q\circ G_{(p,q)}(x) = pQ(x)\bar{q}.$$
The surjective map
$\begin{array}{ccl}
S^3_{\mathbb{H}}\times S^3_{\mathbb{H}} & \longrightarrow & SO(4)\\
(p,q) & \longmapsto & G_{(p,q)}
\end{array}$
is a 2-sheeted covering map (since $G_{(-p,-q)}= G_{(p,q)}$). Explicitely
we have, 

$$G_{(p,q)}x = L_pR_{\bar{q}}x = R_{\bar{q}}L_px,$$
where, denoting $p=p^1+ip^2+jp^3+kp^4$ and $q=q^1+iq^2+jq^3+kq^4$,

$$L_p = p^1\1_4 + p^2L_i+p^3L_j+p^4L_k$$
is the left multiplication by $p$ in $\mathbb{H}$,

$$R_{\bar{q}} = q^1\1_4 -q^2R_i-q^3R_j-q^4R_k$$
is the right multiplication by $\bar{q}$ in $\mathbb{H}$ (notice that 
$R_{\bar{q}}R_{\bar{q'}} = R_{\overline{qq'}}$),
and

$$
L_i = \left( \begin{array}{cccc}
0 & -1 & 0 & 0\\
1 & 0 & 0 & 0\\
0 & 0 & 0 & -1\\
0 & 0 & 1 & 0
\end{array} \right) ,\ 
L_j = \left( \begin{array}{cccc}
0 & 0 & -1 & 0\\
0 & 0 & 0 & 1\\
1 & 0 & 0 & 0\\
0 & -1 & 0 & 0
\end{array} \right) ,\
L_k = \left( \begin{array}{cccc}
0 & 0 & 0 & -1\\
0 & 0 & -1 & 0\\
0 & 1 & 0 & 0\\
1 & 0 & 0 & 0
\end{array} \right) ,
$$
$$
R_i = \left( \begin{array}{cccc}
0 & -1 & 0 & 0\\
1 & 0 & 0 & 0\\
0 & 0 & 0 & 1\\
0 & 0 & -1 & 0
\end{array} \right) ,\ 
R_j = \left( \begin{array}{cccc}
0 & 0 & -1 & 0\\
0 & 0 & 0 & -1\\
1 & 0 & 0 & 0\\
0 & 1 & 0 & 0
\end{array} \right) ,\
R_k = \left( \begin{array}{cccc}
0 & 0 & 0 & -1\\
0 & 0 & 1 & 0\\
0 & -1 & 0 & 0\\
1 & 0 & 0 & 0
\end{array} \right) .
$$
Then, from $\vec{\cal G}\simeq \{ G\in SO(4)/[G,L_i]=0\}$, we obtain

$$\vec{\cal G} \simeq {\cal G}_0. {\cal G}_2,$$
where ${\cal G}_0=\{R_{\bar{q}} = q^1\1_4- q^2R_i-q^3R_j-q^4R_k/q\in 
S^3_{\mathbb{H}}\}$ and
${\cal G}_2=\{L_p = p^1\1_4+ p^2L_i/p\in S^1_{\C}\}$. Notice that, 
for any
$G\in \vec{\cal G}$, there exists $(G_0, G_2) \in {\cal G}_0\times {\cal G}_2$, 
such that $G=G_0 G_2=G_2 G_0$, and $(G_0,G_2)$ is unique up to change of sign.
\lf

Alternatively, we can describe $\vec{\cal G}$ using the isomorphism
\[
\begin{array}{cccl}
C: & \vec{E}^4\simeq \R[4] & \longrightarrow & \C^2\\
 & x^1\epsilon _1+x^2\epsilon _2+x^3\epsilon _3+x^4\epsilon _4\simeq ( 
x^1, x^2, x^3, x^4) & \longmapsto &  (x^1+ ix^2, x^3+ix^4),
\end{array}
\]
which is holomorphic from $(\vec{E}^4, J)$ to $\C^2$. Through
that identification, $\vec{\cal G}$ corresponds to $U(2)$, ${\cal G}_0$ to 
$SU(2)$ and ${\cal G}_2$ to
$$\left\{\left( \begin{array}{cc}
e^{i\theta \over 2} & 0\\
0 & e^{i\theta \over 2}
\end{array} \right) / \theta \in \mathbb{R}\right\}\simeq U(1).$$

It is useful to keep in mind these representations. However, we shall
mostly represent $\vec{\cal G}$ as a subgroup of the $4\times 4$
matrices ${\cal M}(4,\mathbb{R})$ (which we can also identify with a subgroup 
of
${\cal M}(5,\mathbb{R})$, see below), since several complex structures will be
involved.\\

The group ${\cal G}$ is just the semidirect product $\vec{\cal G}\ltimes 
\R[4]$. If $G,G'\in \vec{\cal G}$ and $T,T'\in \R[4]$, the product 
is
$(G,T).(G',T') = (GG',GT'+T)$. This group is embedded in ${\cal 
M}(5,\mathbb{R})$ through

$$(G,T) \longmapsto
\left( \begin{array}{cc}
G & T\\0 & 1
\end{array} \right).$$
We shall call $G$ the {\em rotation component} of $(G,T)$ and $T$ the {\em
translation component} of $(G,T)$. (Notice that we also have the
representation ${\cal G}\simeq U(2)\ltimes \C^2$.)\\

%\subsection{Lie algebras}

The Lie algebra of ${\cal G}$ will be identified with

$$\mathfrak{g} =\left\{ (aL_i+b^1R_i+b^2R_j+b^3R_k, t)/a,b^1,b^2,b^3\in 
\mathbb{R}, t\in \R[4]\right\}.$$
The Lie bracket of two elements $(\eta ,t)$, $(\tilde{\eta}, \tilde{t})\in 
\mathfrak{g}$ is
\[
[ (\eta ,t), (\tilde{\eta}, \tilde{t}) ] = (\eta \tilde{\eta} 
-\tilde{\eta}\eta  , \eta  \tilde{t} - \tilde{\eta}t ).
\]
We denote by $\mathfrak{g}_0$ the Lie algebra of ${\cal G}_0$,
generated by $(R_i,0)$, $(R_j,0)$ and $(R_k,0)$,
and $\mathfrak{g}_2$ the Lie algebra of ${\cal G}_2$, generated by $(L_i, 0)$.
Then the Lie algebra of $\vec{\cal G}$ is $\vec{\g} = \g _0\oplus \g _2$.

\subsection{Action on the Lagrangian Stiefel manifold and on the 
Lagrangian Grassmannian}	\label{Stiefel}

Let us define the Lagrangian Grassmannian $\textit{Gr}_{lag}$ to be the 
set of all oriented 2-dimensional Lagrangian subspaces of $\vec{E}^4$, and
the Lagrangian Stiefel manifold by
\[
\textit{Stief}_\mathit{lag}=\{(e_1,e_3) \in \vec{E}^4\times \vec{E}^4/ 
|e_1|=|e_3|=1, \langle e_1,e_3\rangle =0, \omega (e_1,e_3)=0\}.
\]
Notice that $\textit{Stief}_\mathit{lag}$ is nothing but the set of oriented 
orthonormal bases of planes in $\textit{Gr}_{lag}$.
Actually, we may identify $\textit{Stief}_\mathit{lag}$ with ${\cal F}$ by the 
following: to each basis $(e_1,e_2,e_3,e_4)\in {\cal F}$, we associate
$(e_1,e_3)$ in $\textit{Stief}_\mathit{lag}$. Conversely, we associate to each  
$(e_1,e_3)\in \textit{Stief}_\mathit{lag}$ the frame $(e_1,e_2,e_3,e_4)$ such that
$e_2=L_ie_1$ and $e_4=L_ie_3$. (Through the identification 
$\vec{E}^4\simeq \C^2$, it just amounts to say that $(e_1,e_3)$ is a 
Hermitian basis
of $\C^2$ over $\C$ if and only if $(e_1,ie_1,e_3,ie_3)$ is an 
orthonormal basis of $\C^2$ over $\mathbb{R}$.)

Now the group $\vec{\cal G}$ acts freely and transitively on ${\cal F}$,
i. e., for any $(e_1,e_2,e_3,e_4)\in {\cal F}$, there exists a unique
$G\in \vec{\cal G}$ such that $(e_1,e_2,e_3,e_4) = (\epsilon _1, \epsilon 
_2, \epsilon _3, \epsilon _4)G$. To prove that, it
suffices to realize that the columns of $G$ are  the
components of each vector $e_i$ in the basis $(\epsilon _1, \epsilon _2, 
\epsilon _3, \epsilon _4)$. Hence
$\vec{\cal G}$ acts freely and transitively on $\textit{Stief}_\mathit{lag}$ as 
well, and transitively on $\textit{Gr}_{lag}$:
if $(e_1,e_3) \in \textit{Stief}_\mathit{lag}$ and $G\in \vec{\cal G}$ we shall 
denote $(Ge_1,Ge_3)$ its image by $G$.\\

An important object for the study of Hamiltonian stationary surfaces is 
the {\em Lagrangian angle map}
$\Theta :\textit{Stief}_\mathit{lag} \longrightarrow \R/2\pi \Z$. For 
any $(e_1,e_3)\in \textit{Stief}_\mathit{lag}$, let $G$ be the unique
element in $\vec{\cal G}$ such that $G\epsilon _1 = e_1$ and $G\epsilon _3 
= e_3$, i.e.
$(e_1,L_ie_1,e_3,L_ie_3) = (\epsilon _1, \epsilon _2, \epsilon _3, 
\epsilon _4)G$.
Then, viewing $G$ as a matrix in $U(2)$, we may compute its determinant: 
it is a complex number of modulus one, which
we denote $e^{i\Theta (e_1,e_3)}$. It builds up the {\em Lagrangian angle 
map} $\Theta$. Alternatively, we may decompose $G=G_0 G_2$, 
where $G_0 \in {\cal G}_0$ and $G_2 = e^{\Theta (e_1,e_3) L_i\over 2} 
\in {\cal G}_2$\footnote{in particular it proves that 
$\Theta (GK\epsilon_1,GK\epsilon_3) = \Theta (G\epsilon _1,G\epsilon _3)$, 
$\forall G\in \vec{\cal G}$, $\forall K\in {\cal G}_0$.}.
A last definition is given by

$$(dx^1+idx^2)\wedge (dx^3+idx^4)(e_1,e_3) = e^{i\Theta (e_1,e_3)}.$$
One can check easily that $\Theta (e_1,e_3)$ does not change if we replace 
$(e_1,e_3)$ by another direct
orthonormal basis of the oriented Lagrangian plane  spanned by 
$(e_1,e_3)$. Hence it defines a map from
$\textit{Gr}_{lag}$  to $\R/2\pi \Z$ which we shall also denote 
$\Theta$.

Lastly, in the following we shall abuse notations and identify vectors
$x= x^1\epsilon _1+x^2\epsilon _2+x^3\epsilon _3+x^4\epsilon _4$ in 
$\vec{E}^4$ with column matrices
$\left( \begin{array}{c}x^1 \\ x^2 \\ x^3 \\ x^4\end{array} \right)$. In 
particular we let
$\epsilon _1 =\left( \begin{array}{c}1 \\ 0 \\ 0 \\ 0\end{array} \right)$,
$\epsilon _2 =\left( \begin{array}{c}0 \\ 1 \\ 0 \\ 0\end{array} \right)$,
$\epsilon _3 =\left( \begin{array}{c}0 \\ 0 \\ 1 \\ 0\end{array} \right)$,
$\epsilon _4 =\left( \begin{array}{c}0 \\ 0 \\ 0 \\ 1\end{array} \right)$,
$\epsilon  = {1\over 2} \left( \begin{array}{c}1 \\ 0 \\ -i \\ 
0\end{array} \right)$ and
$\bar{\epsilon}  = {1\over 2} \left( \begin{array}{c}1 \\ 0 \\ i \\ 
0\end{array} \right)$,
$L_i\epsilon  = {1\over 2} \left( \begin{array}{c}0 \\ 1 \\ 0 \\ 
-i\end{array} \right)$ and
$L_i\bar{\epsilon}  = {1\over 2} \left( \begin{array}{c}0 \\ 1 \\ 0 
\\ i\end{array} \right)$.

\subsection{Moving frames for conformal Lagrangian immersions}

Let us consider a smooth conformal Lagrangian immersion of a simply
connected open domain $\Omega$ of $\C\simeq \mathbb{R}^2$, $X:\Omega 
\longrightarrow E^4$.
We shall denote $z=x+iy\simeq (x,y)$ the coordinates on $\R ^2$.
We let $f:\Omega \longrightarrow \R $, such that $e^{f(z)} = 
|{\partial X\over \partial x}| = |{\partial X\over \partial y}|$ and we set
$e_1(z) = e^{-f(z)} {\partial X\over \partial x}(z)$ and
$e_3(z) = e^{-f(z)} {\partial X\over \partial y}(z)$, so that
\[
dX = e^f (e_1dx + e_3dy ),
\]
and then, $X$ is a conformal Lagrangian immersion if and only if for
$z\in \Omega$, $(e_1(z),e_3(z))$ is in $\textit{Stief}_\mathit{lag}$.
Without loss of generality, we will normalize $X$ by
assuming that $X(z_0) = 0$ and $(e_1(z_0), e_3(z_0)) = (\epsilon 
_1,\epsilon _3)$, for some fixed point $z_0\in \Omega$.
We let ${\cal X}$ to be the set of such conformal Lagrangian immersions.
We denote $e_2(z) := L_ie_1(z)$ and $e_4(z) := L_ie_3(z)$.
Therefore the system $e(z):=(e_1(z), e_2(z), e_3(z), e_4(z))$ belongs to 
${\cal F}$: let  $F_X(z)\in \vec{\cal G}$ such that
$e(z) = (\epsilon _1, \epsilon _2, \epsilon _3, \epsilon _4)F_X(z)$ and 
(abusing notations)
let $X(z)$ be the column vector of the components of $X$ in the basis 
$(\epsilon _1, \epsilon _2, \epsilon _3, \epsilon _4)$.
Then we construct a map $\tilde{X}:\Omega \longrightarrow {\cal G}$ 
lifting $X$, defined by

$$\tilde{X}(z) = \left( \begin{array}{cc} F_X(z) & X(z) \\ 0 & 1 
\end{array} \right) \simeq (F_X(z),X(z)).$$
We shall call $\tilde{X}$ the \emph{fundamental lift} of $X$. According to our 
normalization, we have $\tilde{X}(z_0) = (\1 ,0)$.
The Maurer-Cartan form of $\tilde{X}$ is
\[
\tilde{X}^{-1}d\tilde{X} = ( F_X^{-1}dF_X, F_X^{-1}dX).
\]
It is a 1-form with coefficients in $\g$, with the property that its 
translation component has the form
\begin{equation}\label{trans1}
F_X^{-1}dX =  e^f (\epsilon _1 dx +  \epsilon _3 dy)  = e^f (\epsilon dz 
+  \bar{\epsilon} d\bar{z}).
\end{equation}
The key idea in the following will be to study suitably defined lifts of
conformal Lagrangian immersions instead of immersions themselves - which 
has the effect of decreasing by one the order of
the PDE. One could use the fundamental lift. We shall however enlarge the
possibilities as follows:

\begin{definition}\label{LCLI}
A {\em lifted conformal Lagrangian immersion} (LCLI) is  a map 
$U=(F,X):\Omega \longrightarrow {\cal G}$,
satisfying one of the three following equivalent hypotheses.
\begin{itemize}
\item [a)] $U(z) = (F_X(z),X(z)).(K(z)^{-1},0) = (F_X(z)K(z)^{-1}, X(z))$ 
where $X$ is a conformal Lagrangian
immersion, $(F_X,X)$ is its fundamental lift and $K\in {\cal 
C}^{\infty}_{\star}(\Omega , {\cal G}_0) =
\{K\in {\cal C}^{\infty}(\Omega , {\cal G}_0)/ K(z_0) = \1 \}$.
\item [b)] $U(z_0) = (\1 ,0)$ and the translation component of the 
Maurer-Cartan form $U^{-1}dU$ has the form
\[
t = F^{-1}dX = e^f K(\epsilon _1dx +\epsilon _3dy)  = e^f K(\epsilon dz 
+ \bar{\epsilon}d\bar{z}),
\]
where $K\in {\cal C}^{\infty}(\Omega , {\cal G}_0)$ and $f\in {\cal 
C}^{\infty}(\Omega ,\R )$.
\item [c)] $U(z_0) = (\1 ,0)$ and $X$ is a conformal Lagrangian immersion 
and, $\forall z\in \Omega$,
\[
\Theta \left( {\partial X\over \partial x}(z),{\partial X\over \partial 
y}(z)\right)=
\Theta (F(z)\epsilon _1,F(z)\epsilon _3),
\]
where we let $\Theta \left( {\partial X\over \partial x}(z),{\partial 
X\over \partial y}(z)\right)$ to be the value of $\Theta$ on the
oriented Lagrangian plane spanned by ${\partial X\over \partial x}(z)$ 
and  ${\partial X\over \partial y}(z)$.
\end{itemize}
We shall denote by ${\cal GX}$ the set of all LCLI's.
\end{definition}
{\em Proof of the equivalence between  a), b) and c)}.
{\bf {\em a)} $\Rightarrow$ {\em b)}}: it is a direct computation.
{\bf {\em b)} $\Rightarrow$ {\em c)}}: from {\em b)}, it follows that 
${\partial X\over \partial x} = e^fFK\epsilon _1$
and ${\partial X\over \partial y} = e^fFK\epsilon _3$, therefore, using 
the remark in the footnote, section~\ref{Stiefel},
$\Theta \left( {\partial X\over \partial x}(z),
{\partial X\over \partial y}(z)\right)=
\Theta (FK\epsilon _1, FK\epsilon _3) = 
\Theta (F\epsilon _1, F\epsilon _3)$.
{\bf {\em c)} $\Rightarrow$ {\em a)}}: let $(F_X,X)$ be the fundamental 
lift of $X$ and let
$K = F^{-1}F_X\in {\cal C}^{\infty}(\Omega , \vec{\cal G})$.
Then a computation shows that the relation
$\Theta \left( {\partial X\over \partial x}(z),{\partial X\over \partial 
y}(z)\right)=\Theta (F\epsilon _1, F\epsilon _3)$
is equivalent to $\Theta (F_X\epsilon _1, F_X\epsilon _3) = \Theta 
(F_XK^{-1}\epsilon _1, F_XK^{-1}\epsilon _3)$ and using the remark in
the footnote of section 2.3, this implies that $K\in {\cal 
C}^{\infty}_{\star}(\Omega ,{\cal G}_0)$. \bbox

For any simply connected domain $\Omega$ and for any conformal Lagrangian 
immersion $X:\Omega \longrightarrow E^4$,
we shall lift the Lagrangian angle map and define a map $\beta :\Omega 
\longrightarrow \R $, such that
$\forall z\in \Omega$, $\Theta \left( {\partial X\over \partial 
x}(z),{\partial X\over \partial y}(z)\right) =\beta (z)$ modulo
$2\pi$. The (lifted) Lagrangian angle map $\beta$ of a LCLI $U$ is 
characterised by the - unique up to sign - decomposition
$U(z) = (e^{\beta (z)L_i/2} M_0(z),X(z))$, where 
$M_0 \in {\cal C}^{\infty}(\Omega ,{\cal G}_0)$. In the following,
for any $X\in {\cal X}$, we shall choose $\beta$ to be the unique 
Lagrangian angle map such that $\beta (z_0) = 0$.

\begin{remark} \em
It is clear that the gauge group ${\cal C}^{\infty}_{\star}(\Omega ,{\cal 
G}_0)$ acts on the right on ${\cal GX}$ and that
the quotient of ${\cal GX}$ under this gauge action coincides with ${\cal 
X}$. Furthermore, in a given gauge orbit, there are three special
LCLI's: the fundamental lift $\tilde{X}$ and a pair of lifts such that 
$F\in {\cal C}^{\infty}(\Omega ,{\cal G}_2)$,
namely
\[
U_+(z) = (e^{\beta L_i/2},X)\hbox{ and }U_-(z) = (-e^{\beta L_i/2},X).
\]
(Note that by a change $\beta \rightarrow \beta +2\pi$ of the choice of 
the determination of $\beta$, $U_+$ and $U_-$
would be exchanged.) We call $U_+$ and $U_-$ the {\em spinor lifts}.
\end{remark}

\noindent We have the following characterization of Hamiltonian stationary 
surfaces (see \cite{O1}, \cite{ScW}).

\begin{theorem}\label{harmonic}
Let $X:\Omega \longrightarrow E^4$ be a conformal Lagrangian immersion, 
then, $X$ is Hamiltonian stationary if and only
if the Lagrangian angle map is a harmonic function on the surface image, 
or equivalently
\[
\Delta \beta = 0\hbox{ on  }\Omega .
\]
\end{theorem}
Thus we are lead to study LCLI's with harmonic Lagrangian angle map.
Alternatively, we can isolate the differential $d\beta$ by a decomposition 
of the Maurer-Cartan form of
$U=(e^{\beta L_i/2} M_0,X)$ (according to $\mathfrak{g} = 
\mathfrak{g}_2\oplus \mathfrak{g}_0 \oplus (0,\R[4])$),

\begin{equation}\label{MCform}
\alpha = U^{-1}dU = {d\beta \over 2}(L_i,0) + (M_0^{-1} dM_0,0)
+ (0,t).
\end{equation}
Therefore, we may study connection 1-forms $\alpha \in {\cal 
C}^{\infty}(\Omega ,T^{\star}\R ^2\otimes \g )$
on a simply connected domain $\Omega$ which satisfy relation 
(\ref{MCform}) with harmonic $d\beta$, and the zero curvature equation

\begin{equation}\label{curv}
d\alpha +\alpha \wedge \alpha = 0,
\end{equation}
a necessary and sufficient condition for the existence of a map $U:\Omega 
\longrightarrow {\cal G}$ such that
$dU = U.\alpha$; furthermore, $U$ is unique, if we assume also the 
condition

\begin{equation}\label{base}
U(z_0) = (\1 ,0),\hbox{ for some fixed point }z_0\in \Omega .
\end{equation}
We shall concentrate in the following on this last characterization.

\begin{remark} \em
The gauge action of ${\cal C}^{\infty}_{\star}(\Omega ,{\cal G}_0)$ on 
${\cal GX}$ induces an action on Maurer-Cartan 1-forms
described by

$$(\eta ,t)\ \mapsto \ (K\eta K^{-1} -dKK^{-1}, Kt).$$
In any orbit of this gauge action, the fundamental lift 
$\tilde{X}=(F_X,X)=(e^{\beta L_i/2}M_X,X)$ has the Maurer-Cartan form

$$\tilde{X}^{-1}d\tilde{X} = (F_X^{-1}dF_X,0)+(0,e^f(\epsilon dz 
+\bar{\epsilon}d\bar{z})) =
{d\beta \over 2}(L_i,0) + (M_X^{-1}dM_X,0) + (0,e^f(\epsilon dz + 
\bar{\epsilon}d\bar{z})),$$
i. e. with ``simplest'' translation component, whereas the spinor lifts 
has the Maurer-Cartan forms

$$U_{\pm}^{-1}dU_{\pm} = {d\beta \over 2}(L_i,0) + (0,0) + (0,\pm 
e^{-\beta L_i/2}dX),$$
i. e. with zero $\g _0$ component.
\end{remark}

\subsection{Splitting the Lie algebra}

Our aim will be to refine the decomposition given in (\ref{MCform}). We 
introduce the following automorphism
$\tau$ acting on ${\cal G}$ through conjugation by $(-L_j,0)$, i. e. 

$$\tau (G,T) = (-L_j,0) (G,T) (-L_j,0)^{-1} = (-L_jGL_j, -L_jT).$$
It induces a linear action on $\mathfrak{g}$, which diagonalizes on 
$\mathfrak{g}^{\C} = \mathfrak{g}\otimes \C$,
with eigenvalues $i^{-1}$, $i^0$, $i^1$ and $i^2$, since $\tau ^4 = \1$. 
For $k=-1,0,1,2$, we denote by $\mathfrak{g}_k^{\C}$
the eigenspace of $\tau$ for the eigenvalue $i^k$, and we have

\begin{itemize}
\item for the eigenvalue $-i$, $\mathfrak{g}_{-1}^{\C} = (0, 
\C\epsilon \oplus \C L_i\bar{\epsilon})$ (notice that
$\C\epsilon \oplus \C L_i\bar{\epsilon}$ is the 
$(-i)$-eigenspace of $-L_j$),
\item for the eigenvalue 1, $\mathfrak{g}_0^{\C} = 
\mathfrak{g}_0\otimes \C$,
where $\g _0$ is the Lie algebra of ${\cal G}_0$,
\item for the eigenvalue $i$, $\mathfrak{g}_1^{\C} = 
(0,\C\bar{\epsilon}  \oplus \C L_i\epsilon )$ (notice that
$\C\bar{\epsilon}  \oplus \C L_i\epsilon $ is the 
$i$-eigenspace of $-L_j$),
\item for the eigenvalue -1, $\mathfrak{g}_2^{\C} 
=\mathfrak{g}_2\otimes \C$,
where $\g _2$ is the Lie algebra of ${\cal G}_2$.
\end{itemize}

We also have the following characterization of the $\pm i$-eigenspaces.
\begin{lemma}  \label{SU2action}
The group $\R_+^* \times {\cal G}_0$ acts freely and transitively on
the $\pm i$-eigenspaces of $L_j$ minus 0; in particular the $i$-eigenspace 
of $L_j$, $\C\epsilon \oplus \C L_i\bar\epsilon$, coincides with 
the orbit of $\epsilon$ and the $-i$-eigenspace of $L_j$,
$\C\bar\epsilon \oplus \C L_i\epsilon$, coincides with the orbit of 
$\bar\epsilon$.
\end{lemma}
\emph{Proof. }
Since ${\cal G}_0$ commutes with $L_j$,
it preserves its eigenspaces. We now prove the freeness and 
transitivity of the action. Let $\xi = a \epsilon + bL_i \beps
= {1\over 2}(a,b,-ia,ib)$
be an eigenvector associated to the eigenvalue $i$ (for the other
eigenspace use conjugation). If $H \in \R_+^* \times {\cal G}_0$
maps $\epsilon$ to $\xi$, then we infer necessary conditions:
$H\epsilon = {1\over 2} (H\epsilon_1 - iH\epsilon_3) = \xi$;
thus $H\epsilon_1 = {1\over 2}\Re[\xi]$ and $H\epsilon_3 = -{1\over 
2}\Im[\xi]$.
Since we want $H$ to commute with $L_i$, $H\epsilon_2 = H L_i\epsilon_1
= {1\over 2}L_i \Re[\xi] = {1\over 2}\Re[L_i \xi]$, and $H\epsilon_4 = 
-{1\over 2}\Im[L_i \xi]$.
So $H$ is uniquely determined. Check easily that $( H\epsilon _1, 
H\epsilon _2, H\epsilon _3, H\epsilon _4)$
is a conformal basis of \R[4], and write $H=rK$ for some
isometry $K$ and some $r\in \R ^{\star}_+$. By construction $K$ 
commutes with $L_i$, so
$K$ belongs in $\vec{\cal G}$. It cannot
have any nontrivial component in ${\cal G}_2$ otherwise $\xi$ would not be
a eigenvector of $L_j$ (whose eigenspaces are not stable under $L_i$).
Thus $K$ belongs to  ${\cal G}_0$. Therefore there exists a unique $H \in 
\R_+^* \times {\cal G}_0$
sending $\epsilon$ to $\xi$.

Notice that the action of $\R_+^* \times {\cal G}_0$ coincides with the 
right action of $\mathbb{H}^{\star}$ on $\C\epsilon \oplus 
\C L_i\bar{\epsilon}$ as a subset of $\mathbb{H}\otimes \C$. 
\bbox

Using the decomposition
$\mathfrak{g}^{\C} = \mathfrak{g}_{-1}^{\C} \oplus 
\mathfrak{g}_0^{\C} \oplus \mathfrak{g}_1^{\C}
\oplus \mathfrak{g}_2^{\C}$, we define the projection mapping 
$[.]_k:\mathfrak{g}^{\C}\longrightarrow \mathfrak{g}_k^{\C}$.
Then denoting $\alpha _k:=[\alpha ]_k$, we have

\begin{equation}\label{alpha}
\alpha = \alpha _{-1}+\alpha _0+\alpha _1+\alpha _2.
\end{equation}
We now substitute (\ref{alpha}) in (\ref{curv}). We use the relations
$[\mathfrak{g}_k^{\C}, \mathfrak{g}_l^{\C}] \subset 
\mathfrak{g}_{(k+l) \bmod 4}^{\C}$ and
$[\mathfrak{g}_{\pm 1}^{\C}, \mathfrak{g}_{\pm 1}^{\C}] = 
[\g_0,\g_2]=[\g_2,\g_2]=0$. The projection
of the resulting equation on each eigenspace  gives us four relations

\begin{equation}\label{4rel}
\left\{
\begin{array}{ccc}
d\alpha _{-1} + [\alpha _{-1}\wedge \alpha _0] + [\alpha _1\wedge \alpha 
_2] & = 0,\\
d\alpha _0 +  {1\over 2} [\alpha _0\wedge \alpha _0] & = 0,\\
d\alpha _1 + [\alpha _0\wedge \alpha _1] + [\alpha _{-1}\wedge \alpha _2] 
& = 0,\\
d\alpha _2  & = 0,
\end{array} \right.
\end{equation}
where $[\alpha _a\wedge \alpha _b] = \alpha _a\wedge \alpha _b + \alpha 
_b\wedge \alpha _a$. We further
decompose each form $\alpha _k$ as $\alpha _k = \alpha _k' +\alpha _k''$, 
with
$\alpha _k' = \alpha _k({\partial \over \partial z})dz$ and $\alpha _k'' = 
\alpha _k({\partial \over \partial \bar{z}})d\bar{z}$.
We remark that, because $\alpha$ derives from a LCLI,
$\alpha _{-1} +\alpha _1 = (0,e^{f(z)}K(z)(\epsilon dz + 
\bar{\epsilon} d\bar{z}))$ and hence, Lemma \ref{SU2action}
implies that

\begin{equation}\label{alpha-1}
\alpha _{-1} = (0,e^{f(z)}K(z)\epsilon dz) = \alpha _{-1}'\hbox{ and 
}\alpha _{-1}'' = 0,
\end{equation}
and similarly,

\begin{equation}\label{alpha1}
\alpha _1 =\alpha _1''\hbox{ and }\alpha _1' = 0.
\end{equation}
Thus,

$$
\alpha = \alpha _2' + \alpha _{-1}' + \alpha _0 + \alpha _1'' + \alpha 
_2''.
$$
Now we exploit (\ref{alpha-1}) and  (\ref{alpha1}) in (\ref{4rel}) and we 
obtain

$$\left\{ \begin{array}{ccc}
d\alpha _{-1}' + [\alpha _{-1}'\wedge \alpha _0] + [\alpha _1''\wedge 
\alpha _2'] & = & 0,\\
d\alpha _0 +  {1\over 2} [\alpha _0\wedge \alpha _0] & = & 0,\\
d\alpha _1'' + [\alpha _0\wedge \alpha _1''] + [\alpha _{-1}'\wedge \alpha 
_2''] & = & 0,\\
d\alpha _2  & =  & 0.
\end{array} \right. $$
A convenient way to rewrite this system is to introduce a complex 
parameter $\lambda\in \C^{\star}$ and to let

$$\alpha _{\lambda} := \lambda ^{-2} \alpha _2' + \lambda ^{-1} \alpha 
_{-1}' + \alpha _0 + \lambda \alpha _1'' + \lambda ^2\alpha _2'',$$
and then

$$\begin{array}{cccl}
d\alpha _{\lambda} + \alpha _{\lambda}\wedge \alpha _{\lambda} & = & & 
\lambda ^{-2} d\alpha _2'\\
 & & + & \lambda ^{-1} (d\alpha _{-1}' + [\alpha _{-1}'\wedge \alpha _0] + 
[\alpha _1''\wedge \alpha _2'])\\
 & & + & (d\alpha _0  + {1\over 2} [\alpha _0\wedge \alpha _0])\\
 & & + & \lambda (d\alpha _1'' + [\alpha _0\wedge \alpha _1''] + [\alpha 
_{-1}'\wedge \alpha _2''])\\
 & & + & \lambda ^2d\alpha _2''\\
 & = & & \lambda ^{-2} d\alpha _2' + \lambda ^2d\alpha _2''. 
\end{array}$$
We now are in position to prove the 

\begin{theorem}
Assume that $\Omega$ is a simply connected domain of $\C\simeq 
\R ^2$.
Let $\alpha$ be in ${\cal C}^{\infty}(\Omega ,T^{\star}\R ^2\otimes 
\mathfrak{g})$. Then
\begin{itemize}
\item $\alpha$ is the Maurer-Cartan form of a LCLI if and only if $d\alpha 
+\alpha \wedge \alpha =0$, 
$\alpha _{-1}'' = \alpha _1' = 0$  and $\alpha _{-1}'\neq 0$, $\alpha 
_1''\neq 0$,
\item furthermore, it corresponds to some Hamiltonian stationary immersion 
if and only if
the {\em extended Maurer-Cartan form}
$\alpha _{\lambda} = \lambda ^{-2} \alpha _2' + \lambda ^{-1} \alpha 
_{-1}' + \alpha _0 + \lambda \alpha _1'' + \lambda ^2\alpha _2''$
satisfies
\begin{equation}\label{lambda}
d\alpha _{\lambda} + \alpha _{\lambda}\wedge \alpha _{\lambda}= 0,\ 
\forall \lambda \in \C^{\star}.
\end{equation}
\end{itemize}
\end{theorem}
{\em Proof.}  First, according to Definition \ref{LCLI}, {\em b)}, $\alpha$ 
will be the Maurer-Cartan form of a LCLI if
and only if $d\alpha +\alpha \wedge \alpha =0$ and

$$(\alpha _{-1} + \alpha _1)({\partial \over \partial z}) = 
(0,e^{f(z)}K(z)\epsilon )\hbox{ and }
(\alpha _{-1} + \alpha _1)({\partial \over \partial \bar{z}}) = 
(0,e^{f(z)}K(z)\bar{\epsilon} ).$$
But, from Lemma \ref{SU2action}, this is equivalent to 

$$(\alpha _{-1} + \alpha _1)({\partial \over \partial z})\in 
\mathfrak{g}_{-1}^{\C}\setminus \{0\} \hbox{ and }
(\alpha _{-1} + \alpha _1)({\partial \over \partial \bar{z}})\in 
\mathfrak{g}_1^{\C}\setminus \{0\} ,$$
or $\alpha _{-1}'' = \alpha _1' = 0$ and $\alpha _{-1}'\neq 0$, $\alpha 
_1''\neq 0$.

Second, the previous computation shows that

$$d\alpha _{\lambda} + \alpha _{\lambda}\wedge \alpha _{\lambda} =
{1\over 2}{\partial ^2\beta \over \partial z\partial \bar{z}} 
(\lambda ^{-2} - \lambda ^2)(L_i,0)d\bar{z}\wedge dz,$$
which vanishes if and only if the immersion is Hamiltonian stationary, 
according to Theorem~\ref{harmonic}.  \bbox

Notice that it suffices to check Relation (\ref{lambda}) for $\lambda \in 
S^1\subset \C^{\star}$, or even for one value of
$\lambda$ different from $\pm 1$ to ensure 
the Hamiltonian stationary condition.\\

\begin{corollary}
Assume that $\Omega$ is simply connected. Let $\alpha$ be in ${\cal 
C}^{\infty}(\Omega ,T^{\star}\R ^2\otimes \mathfrak{g})$, a 
Maurer-Cartan form
of a Hamiltonian stationary LCLI and $z_0\in \Omega$.
Then for any $\lambda \in S^1$, there exists a unique LCLI $U_{\lambda}\in 
{\cal C}^{\infty}(\Omega ,{\cal G})$
such that

\begin{equation}\label{dU}
dU_{\lambda} = U_{\lambda} \alpha _{\lambda}\hbox{ and }U_{\lambda} (z_0) 
= \1 .
\end{equation}
Thus there is a $S^1$-family of Hamiltonian stationary Lagrangian 
conformal immersions $X_{\lambda}$ given by
$U_{\lambda} = (F_{\lambda}, X_{\lambda})$.
\end{corollary}
{\em Proof.} First the condition $\lambda \in S^1$ ensures that $\alpha 
_{\lambda}$ is $\g$-valued (and not $\g ^{\C}$-valued).
Then equation (\ref{lambda}) is the necessary and sufficient condition for 
the existence of a unique solution to (\ref{dU}). \bbox

Recovering $U_{\lambda}$ from $\alpha _{\lambda}$ can be done in two steps 
(this is
due to the semiproduct structure of ${\cal G}$), namely:
the rotation term can be obtained by solving
$(F_{\lambda}^{-1}dF_{\lambda},0) = [\alpha _{\lambda}]_2 + \alpha_0$; 
recall however that $F_{\lambda}$
is defined only up to gauge transformation, which leaves 
$[\alpha_{\lambda}]_2$
invariant but changes all the other components. The immersion 
$X_{\lambda}$ 
is obtained by solving $(0,F_{\lambda}^{-1}dX_{\lambda}) = \lambda 
^{-1}\alpha_{-1}'+\lambda ^1\alpha_1''$.

The family of solutions $(X_\lambda)_{\lambda \in S^1}$ is quite similar 
to the conjugate family 
of minimal surfaces, also obtained by varying a parameter in $S^1$. 
As in the classical minimal case, this family is in general not 
well-defined if the parameter domain is not simply connected; there
may be period problem (in our setting: non trivial monodromy).
Notice a big difference though: the group involved in the classical 
minimal surface theory is simply \R[3], which unlike \G\/ is commutative.

%%%%%%%%%%%%%%%%%%%%%%%%%%%%%%%%%%%%%%%%%%%%%%%%%%%%%%%%%%%%%%%%%%%%%
%%%%%%%%%%%%%%%%%%%%%%%%%%%%%%%%%% LINEAR PROBLEM  %%%%%%%%%%%%%%%%%%
%%%%%%%%%%%%%%%%%%%%%%%%%%%%%%%%%%%%%%%%%%%%%%%%%%%%%%%%%%%%%%%%%%%%%

\section{An associated linear problem}\label{linearproblem}

In this section we show how a particular choice of gauge (or equivalently
a particular moving frame) reduces the problem to solving successively 
three surprisingly simple linear PDEs, the two first involving 
the conformal structure, the third being simply the integration procedure 
from the connection 1-form to the immersion. We are then in the position 
to describe explicitely all weakly conformal Hamiltonian stationary 
Lagrangian tori.

\subsection{Using the spinor lift}

In this section we still assume that the immersion $X$ is defined on 
some simply connected domain $\Omega$. Hence there is no problem 
in considering a spinor lift $U = (e^{\beta L_i/2},X)$, whose
Maurer-Cartan form is 
\[
\alpha = U^{-1}dU = \left( \frac{1}{2} \pa{\beta}{z} (L_i,0) 
+ (0,u) \right) dz + \left( (0,\bar{u}) + \frac{1}{2} \pa{\beta}{\bz} 
(L_i,0) 
\right) d\bz
\]
where $u = a \epsilon + b L_i \bar{\epsilon}$ for some smooth complex 
valued functions $a,b$ (recall Lemma 1). Equation~(\ref{lambda}) yields 
the condition
\begin{equation}	\label{lin1}
\Delta \beta = 0 .
\end{equation}
Using the fact that in our connexion form $\alpha_0 = 0$,
the only other condition in~(\ref{lambda}) is another linear PDE:
\[
\left( 0,\pa{u}{\bz} \right) + [(0,\bar{u}), \frac{1}{2} 
\pa{\beta}{z}(L_i,0)] = 0;
\]
written more simply:
\begin{equation}	\label{lin2}
\pa{u}{\bz} = \frac{1}{2} \pa{\beta}{z} L_i \bar{u} .
\end{equation}
Finally, once $\beta$ and $u$ are found, $X$ is obtained by 
integrating 
\begin{equation}	\label{lin3}
dX = e^{\beta L_i/2}(u dz + \bar{u}d\bz).
\end{equation}
Notice that the set of solutions of~(\ref{lin2}) is a real vector space;
thus the set of solutions for $X$ is the orbit under \G\/ of
a vector space. Beware also that solving~(\ref{lin2}) does not
garantee that $u$ (and hence the induced metric) will never vanish;
so we may actually obtain weakly conformal solutions. Therefore the 
conclusion:

\begin{theorem}
The Hamiltonian stationary conformal Lagrangian immersions from
a simply connected domain $\Omega$ into $E^4$ are given by solving 
successively three linear partial differential equations
(\ref{lin1}), (\ref{lin2}) and (\ref{lin3}). Then for given 
conformal structure and Lagrangian angle map $\beta$, 
the set of weakly conformal solutions is the \G-orbit of a vector space.
\end{theorem}

\subsection{Hamiltonian stationary Lagrangian tori}

We specialize to the case of Hamiltonian stationary Lagrangian tori.
Let us fix some notations: $\Gamma$ is a lattice in \C,
with dual lattice $\Gamma^* = \{ \gamma \in \C, \< \gamma, \Gamma \> 
\subset \Z \}$ (here $\<.,.\>$ is the usual dot product in 
$\C \simeq \R[2]$); then any torus \T\/ is conformally equivalent to some
$\C/\Gamma$. 
We want to classify Hamiltonian stationary conformal Lagrangian maps $X$
from $\T$ to \R[4], or equivalently their $\Gamma$-periodic lift to the
universal cover, that we will also write abusively $X : \C \to \R[4]$.
This amounts to finding solutions of Equations (\ref{lin1}), (\ref{lin2}) and
(\ref{lin3}) which give $\Gamma$-periodic maps on \C.

First one should notice
that the rotation part in the spinor lift $U = (e^{\beta L_i/2},X)$,  
is not $\Gamma$-periodic but only $2\Gamma$-periodic a priori\footnote{
recall that the dual lattice to $2\Gamma$ is just $\frac{1}{2}\Gamma^*$.}; 
so is the translation part of the corresponding Maurer-Cartan form 
and in particular the complex vector $u=e^{-\beta L_i/2}\pa{X}{z}$.
So we will distinguish the ``truly periodic'' solutions, for which 
$e^{\beta L_i/2}$ is $\Gamma$-periodic, from the ``anti-periodic'' ones; 
still, given an anti-periodic solution, its fourfold\footnote{i.e. twofold 
in each lattice direction.} cover is truly periodic.
This detail will become important when we restrict the solutions 
obtained on the universal cover \C\/ to the torus \T. \lf

The solutions of equation~(\ref{lin1}) are particularly simple: since 
$e^{i\beta}$ is periodic,
i.e.  $\beta(z+\Gamma) \equiv \beta(z) \bmod 2\pi$, we have 

\begin{equation}\label{beta}
\beta(z) = 2\pi \<\beta_0,z-z_0\>
\end{equation}
for some  $z_0 \in \Omega$ and $\beta _0\in \Gamma ^*$.
Up to a translation in $z$ we may suppose 
that $z_0=0$.
We see that $e^{\beta L_i/2}$ is $\Gamma$-periodic if and only if 
$\beta_0/2$
belongs to $\Gamma^*$; otherwise $e^{\beta L_i/2}$ is just 
anti-periodic (we will give examples of both cases). 
Now, setting $u = a \epsilon + b L_i \bar{\epsilon}$, 
equation~(\ref{lin2}) is equivalent to
\begin{equation}\label{ab}
\pa{a}{\bz} = - \frac{\pi\bar{\beta}_0}{2} \bar{b} 
\quad \textrm{ and } \quad
\pa{b}{\bz} = \frac{\pi\bar{\beta}_0}{2} \bar{a}.
\end{equation}
A necessary condition for $(a,b)$ to be a solution of (\ref{ab}), is that 
$a$ and $b$ solve the eigenvalue problem:

\begin{equation}\label{egv}
\Delta \psi + \pi^2 |\beta_0|^2 \psi = 0.
\end{equation}
Since $a$ is $2\Gamma$-periodic, it has the
Fourier expansion $a = \sum_{\gamma \in \Gamma^*/2} \hat{a}_{\gamma} 
e^{2i\pi \<\gamma,z\>}$,
and $a$ is a solution of (\ref{egv}) if and only if
all coefficient $\hat{a}_\gamma$ vanish unless $|\gamma| = |\beta_0/2|$.
Notice that except for $\pm \beta_0/2$, existence of such lattice points 
is far from obvious and depends strongly on the conformal structure 
together with the choice of the lattice point $\beta_0$.
We remark now that if $a$ is a solution of (\ref{egv}) and if $b$ is given
by the first equation of (\ref{ab}), then $a$ and $b$ are automatically 
solutions
of the second equation in (\ref{ab}). We deduce the following
(all sums being taken with $\gamma \in  \frac{1}{2}\Gamma^*$)
\[
a(z) = \sum_{|\gamma|=|\beta_0|/2} \hat{a}_{\gamma} e^{2i\pi \<\gamma,z\>} 
\; , \qquad
b(z) = \frac{2}{i\beta_0} \sum_{|\gamma|=|\beta_0|/2} \bar{\gamma}
\overline{\hat{a}_{-\gamma}}  e^{2i\pi \<\gamma,z\>}
= -\frac{2}{i\beta_0} \sum_{|\gamma|=|\beta_0|/2}  \bar{\gamma}
\overline{\hat{a}_{\gamma}}  e^{-2i\pi \<\gamma,z\>}
\]
We conclude that, for $\beta$ given by (\ref{beta}), any solution to 
(\ref{lin2}) has the form $u = \sum _{|\gamma|=|\beta_0|/2} u_{\gamma}$,
where each
$u_{\gamma} = \hat{a}_{\gamma} e^{2i\pi \<\gamma,z\>}\epsilon +
{2i\bar{\gamma}\over \beta _0}\overline{\hat{a}_{\gamma}}  e^{-2i\pi 
\<\gamma,z\>}L_i\bar{\epsilon}$. In other words, the
set of solutions to (\ref{lin2}) is a finite real vector space with basis 
vectors

$$v_{\gamma} = e^{2i\pi \<\gamma,z\>}\epsilon + {2i\bar{\gamma}\over \beta 
_0} e^{-2i\pi \<\gamma,z\>}L_i\bar{\epsilon}
\hbox{ and }w_{\gamma} = ie^{2i\pi \<\gamma,z\>}\epsilon+ 
{2\bar{\gamma}\over \beta _0} e^{-2i\pi \<\gamma,z\>}L_i\bar{\epsilon},$$
for $\gamma \in \{\gamma \in \frac{1}{2}\Gamma^*/ |\gamma|=|\beta_0|/2\}$.

The last step is finding $X$, 
by integrating (\ref{lin3}). Again, assuming that $X(0)=0$,
the set of solutions is a vector space, with the basis vectors
\[
A_{\gamma}(z) = \int _0^z e^{\pi \< \beta_0, \xi \>L_i}(v_{\gamma}(\xi 
)d\xi + \overline{v_{\gamma}(\xi )}d\bar{\xi})
\hbox{ and }
B_{\gamma}(z) = \int _0^z e^{\pi \< \beta_0, \xi \>L_i}(w_{\gamma}(\xi 
)d\xi + \overline{w_{\gamma}(\xi )}d\overline{\xi}),
\]
and any solution has the form
\[
X = \sum_{|\gamma|=|\beta_0|/2} \Re(\hat{a}_{\gamma})A_{\gamma} + 
\Im(\hat{a}_{\gamma})B_{\gamma}.
\]
In the computation of $A_{\gamma}$ and  $B_{\gamma}$, two cases occur: 
either $\gamma  = \pm \frac{1}{2} \beta_0$, and then
$A_{\gamma}$ and $B_{\gamma}$ cannot be periodic, but only pseudo-periodic 
(and both periods cannot compensate).
If  $\gamma  \neq \pm \frac{1}{2} \beta_0$, then $A_{\gamma}$ and 
$B_{\gamma}$ have frequencies
$\gamma \pm \frac{1}{2} \beta_0$ and more precisely:
\[
A_{\gamma} = {4e^{\pi \<\beta _0,z\>L_i}\over \pi} \hbox{Re}\left[ {e^{-2i\pi\<\gamma,z\>}\over \beta _0^2-4\gamma ^2}
\left(  \begin{array}{c}
	\ds -i\gamma \\ \ds -\beta _0/2 \\ 
	\ds  \gamma \\ \ds  -i\beta _0/2
\end{array} \right) \right] ,
\]
and
\[
B_{\gamma} = {4e^{\pi \<\beta _0,z\>L_i}\over \pi} \hbox{Im}\left[  {e^{-2i\pi\<\gamma,z\>}\over \beta _0^2-4\gamma ^2}
\left(  \begin{array}{c}
	\ds -i\gamma \\ \ds -\beta _0/2 \\ 
	\ds  \gamma \\ \ds  -i \beta _0/2
\end{array} \right) \right] .
\]
A necessary and sufficient condition for $X = \sum_{|\gamma|=|\beta_0|/2} 
\Re(\hat{a}_{\gamma})A_{\gamma} + \Im(\hat{a}_{\gamma})B_{\gamma}$ 
to be $\Gamma$-periodic is obviously $\gamma-\frac{1}{2}\beta_0 \in 
\Gamma^*$
(then $\gamma+\frac{1}{2}\beta_0$, $-\gamma+\frac{1}{2}\beta_0$, 
$-\gamma-\frac{1}{2}\beta_0$ automatically belong to $\Gamma^*$).
So we define the set 
\[
\Gamma^*_{\beta_0} = \left\{ \gamma \in \frac{\beta_0}{2}+\Gamma^* 
\textrm{ such that } |\gamma|^2 = \left|\frac{\beta_0}{2}\right|^2 
\textrm{ and } \gamma^2 \neq \left(\frac{\beta_0}{2}\right)^2   \right\}. 
\]

\begin{remark} \em
In the truly periodic case, $\frac{1}{2}\beta_0$ belongs to
$\Gamma^*$, thus $\Gamma^*_{\beta_0}$ is just the intersection of 
the dual lattice with the circle through $\frac{1}{2}\beta_0$,
minus $\pm \frac{1}{2}\beta_0$. 
\end{remark}

\begin{remark} \em 		\label{multG2}
As noted above, multiplication of the solution
by a constant matrix in $\G_2$ is equivalent to a translation in
$z$-space; such a change of variable in turn amounts to multiplying
each $\hat{a}_\gamma$ by a constant (depending on $\gamma$). Furthermore, 
the group action 
of $\R_+^* \times \G_0$ on $X$ descends to a free action on the couples
$(\hat{a}_\gamma,\hat{a}_{-\gamma})$ (transitive on each couple), so that 
all solutions are obtained once for each choice of the $\hat{a}_\gamma$,
up to the obvious $\beta(0) = 0$ assumption.
\end{remark}
We may now conclude by the following classification theorem:
\begin{theorem}
The Hamiltonian stationary weakly conformal Lagrangian immersions from
$\C/\Gamma$ into $E^4$ are characterized by their Lagrangian angle $\beta$
in as much as $\beta_0 = \frac{1}{\pi} \partial \beta/\partial \bz$
belongs to the dual lattice $\Gamma^*$. The set of solutions for a chosen 
$\beta_0$ is the orbit under $\G_2$ of the vector space generated by the 
$A_\gamma$, $B_\gamma$ and translations in 4-space, as $\gamma$ ranges 
over the (possibly empty) set $\Gamma^*_{\beta_0}$. Its dimension -- 
if not empty -- is $2\mathrm{Card}(\Gamma^*_{\beta_0}) + 5$, or 
$2\mathrm{Card}(\Gamma^*_{\beta_0}) -3$ if one identifies solutions in
the same \G-orbit.
\end{theorem}

\begin{remark} \em
It also often happens that a solution 
constructed that way is actually a multiple cover
of another solution with a potentially different conformal type;
indeed the relevant dual lattice is obtained as the one generated by all
$\gamma - \frac{1}{2}\beta_0$, $\gamma + \frac{1}{2}\beta_0$ for
$\gamma \in \Gamma^*_{\beta_0}$ and $\hat{a}_\gamma \neq 0$. As will 
be shown below, truly periodic examples are not always covers of 
antiperiodic ones, though that is that case for square tori:
\end{remark}
\begin{proposition}
Let $\T=\C/\Gamma$ be a square torus; then a truly periodic solution 
$X$ is always a twofold cover of some simpler solution.
\end{proposition}
\emph{Proof. }
For simplicity assume $\Gamma=\Gamma^*$ is just $\Z \oplus \Z$. By hypothesis
$\beta_0$ belongs to $2\Gamma^*$ . We denote $\Delta^*$, dual to $\Delta$, 
the lattice generated by the $\gamma-\frac{1}{2}\beta_0$ and 
$\gamma+\frac{1}{2}\beta_0$ for all relevant\footnote{that is 
$\hat{a}_\gamma \neq 0$.} $\gamma \in \Gamma^*_{\beta_0}$. 
We claim that $\Delta^*$ is a subgroup of $\Delta^*_0 = \{ (n,m) \in 
\Gamma;
n+m \equiv 0 \bmod 2 \} = (1-i)\Z \oplus (1+i)\Z$. Indeed let 
$\gamma = (p,q)$ and $\frac{1}{2}\beta_0 = (p_0,q_0)$ -- both in 
$\Gamma$ --, then $p^2+q^2 = p_0^2 + q_0^2$ implies
\[
p-p_0 \equiv (p-p_0)(p+p_0) \equiv (q-q_0)(q+q_0) \equiv  q-q_0 \bmod 2
\]
Now $\Delta^* \subset \Delta_0^*$ is equivalent to $\Delta_0 
= \frac{1-i}{2}\Z \oplus \frac{1+i}{2}\Z \subset \Delta$, 
so $X$ is $\Delta_0$-periodic, and $\C/\Gamma$ is a double cover
of $\C/\Delta_0$. Notice that $\Delta_0$ is again square. \bbox

\begin{remark} \em
Examples with arbitrarily many frequencies can be 
constructed by taking for instance the square lattice $\Z\oplus \Z$
and $\beta_0 = 2 (p+iq)^{n}$ where $p,q$ are integers. Then
the set $\Gamma^*_{\beta_0}$ contains all the $\pm (p+iq)^{n-k}(p-iq)^k$
for $k$ ranging from $1$ to $n$.
\end{remark}

\subsection{Some toric examples}

\paragraph{A truly periodic example on a non rectangular torus. }
\lf
Set $\omega = e^{i\pi/3}$, and $\Gamma^* = \Z \oplus \omega\Z$.
Let $X$ be a solution with $\beta_0 = 2 \in 2\Gamma^*$ and
non zero coefficients $\hat{a}_{\omega}$, $\hat{a}_{\omega^2}$.
Let $\Delta$ be the lattice of periods of $X$, and $\Delta^*$ its dual.
Then $\Delta^*$ contains $1 = (\omega-\frac{1}{2}\beta_0)
-(\omega^2-\frac{1}{2}\beta_0)$ and $\omega = 
(\omega+\frac{1}{2}\beta_0)-1$.
So $X$ is not a cover of some simpler example. Taking for instance
$\hat{a}_\omega = \hat{a}_{\omega^2} = 1$, we have 
$X(z) = A_\omega(z) + A_{\omega^2}(z)$:
\begin{eqnarray*}
X(z) &=& \frac{e^{2\pi x L_i}}{\pi\sqrt{3}} \left( \begin{array}{c}
	\cos(\pi(x+y\sqrt{3}))  \\
	-\sin(\pi(x+y\sqrt{3}+\frac{1}{3}))  \\
	\sin(\pi(x+y\sqrt{3}))  \\
	\cos(\pi(x+y\sqrt{3}+\frac{1}{3}))
\end{array} \right)
+ \frac{e^{2\pi x L_i}}{\pi\sqrt{3}} \left( \begin{array}{c}
	\cos(\pi(-x+y\sqrt{3}))  \\
	\sin(\pi(-x+y\sqrt{3}-\frac{1}{3}))  \\
	\sin(\pi(-x+y\sqrt{3}))  \\
	-\cos(\pi(-x+y\sqrt{3}-\frac{1}{3}))
\end{array} \right)
\\
 &=& \frac{2}{\pi\sqrt{3}} \left( \begin{array}{c}
	\cos(\pi y\sqrt{3}) \left( 
	\cos(2\pi x)\cos(\pi x)+\sin(2\pi x) \sin(\pi x+\frac{\pi}{3}) \right) 
	\\
	\cos(\pi y\sqrt{3}) \left( 
	\sin(2\pi x)\cos(\pi x)-\cos(2\pi x) \sin(\pi x+\frac{\pi}{3}) \right)  
	\\
	\sin(\pi y\sqrt{3}) \left( 
	\cos(2\pi x)\cos(\pi x)+\sin(2\pi x) \sin(\pi x+\frac{\pi}{3}) \right)  
	\\
	\sin(\pi y\sqrt{3}) \left( 
	\sin(2\pi x)\cos(\pi x)-\cos(2\pi x) \sin(\pi x+\frac{\pi}{3}) \right)
\end{array} \right) .
\end{eqnarray*}
\begin{figure}[h]
\begin{center}
\epsfig{file=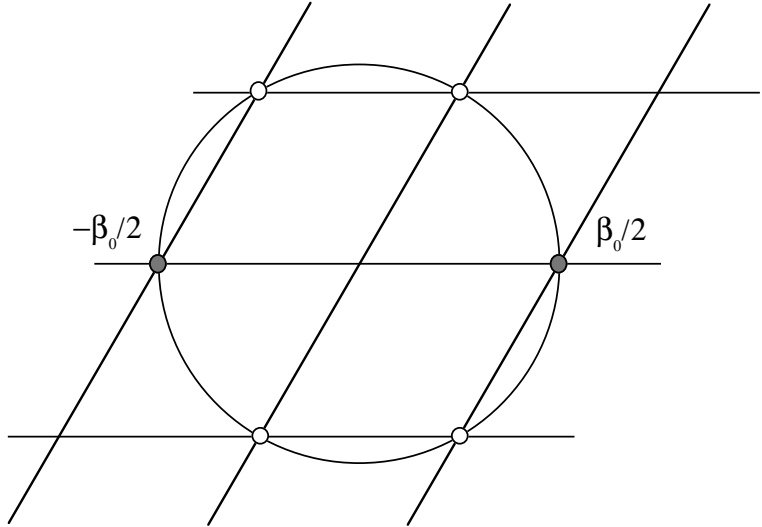, width=10cm}
\end{center}
\caption{A rhombic torus}
\end{figure}

\paragraph{The standard torus and its rectangular counterparts. }
\lf
The simplest -- and until recently (cf~\cite{CU}) -- only known
tori were the product of circles $S^1 \times S^1 \subset \C \times \C
= \R[4]$; more precisely define on the rectangular torus 
$\T = \C/\omega_1 \Z \oplus i\omega_2 \Z$ ($\omega_1,\omega_2 \in \R$)
\[
X(x+iy) = \left( \begin{array}{c}
	\omega_1 \sin (2\pi x/\omega_1)  \\
	-\omega_1 \cos (2\pi x/\omega_1)  \\
	\omega_2 \sin (2\pi y/\omega_2)  \\
	-\omega_2 \cos (2\pi y/\omega_2)
\end{array}  \right)
\]
When $\omega_1=\omega_2=1$, $X(\T)$ is the square torus.

The Lagrangian angle of $X$ is $\beta(z) = 2\pi(x/\omega_1+y/\omega_2)$
so $\beta_0 = \omega_1^{-1} + i \omega_2^{-1}$ which belongs to
$\Gamma^* = \omega_1^{-1}\Z \oplus i \omega_2^{-1} \Z$ but not to
$\frac{1}{2} \Gamma^*$, so we are in the antiperiodic case. Then
$\Gamma^*_{\beta_0} = \{ \frac{1}{2} \bar{\beta}_{0}, 
-\frac{1}{2} \bar{\beta}_{0} \}${~}\footnote{except for some particular 
lattices: if $1+\omega_k^2/\omega_{\ell}^2 = m^2$ for some
$m \in \Z$ then $\pm m \omega_k^{-1}$ also belongs to 
$\Gamma^*_{\beta_0}$.}.
That torus corresponds exactly to $\hat{a}_{\bar{\beta}_0/2} 
= \hat{a}_{-\bar{\beta}_0/2} = \pi$). As noted in remark~\ref{multG2}, 
other choices of those coefficients amount to multiplying $X$ by an element 
in 
$\R_+^* \times \G_0$.
\begin{figure}[h]
\begin{center}
    \epsfig{file=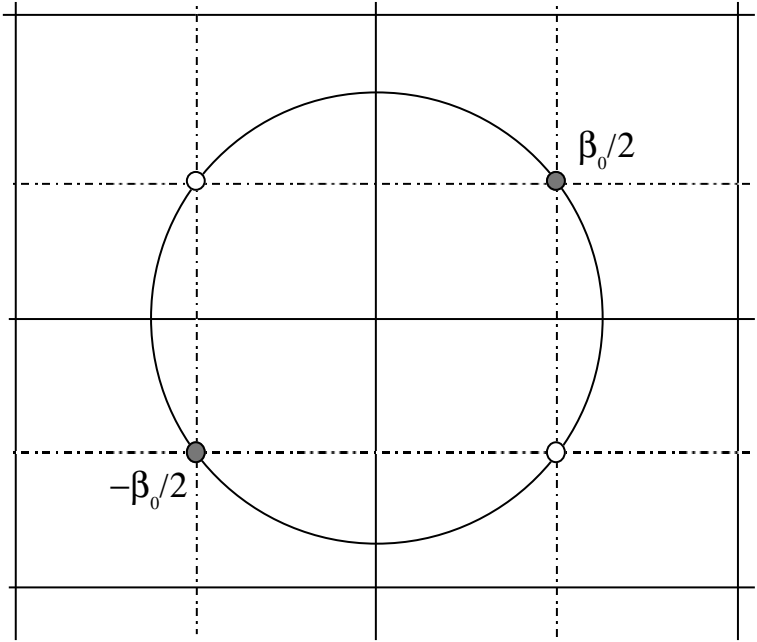, width=10cm}
\end{center}
\caption{The standard torus}
\end{figure}

\paragraph{The examples of I. Castro and F. Urbano.}
\lf
In a recent article a new 3-parameter family of Hamiltonian stationary
Lagrangian surfaces was described, some of them giving rise to tori
(when the periodicity conditions were satisfied). The construction was 
based on geometric properties of parallel lifts on the 3-sphere,
instead of the more analytic methods used here. These examples also
satisfied the rigidity property of being invariant under a 
one-parameter group of isometries, which characterizes them among
tori without parallel mean curvature vector of course (see the 
examples above). Though we will not describe explicitely the examples
(we refer the interested Reader to~\cite{CU}), we
indicate how they fit in our classification and how their 
properties are linked with a special lattice structure. From now on 
we will use their notations.

Let  $X_{\theta,\beta}^\alpha$ be an immersion with real parameters
$\alpha,\beta,\theta \in [0,\pi/2) \times (0,\pi/2) \times 
(-\pi/2,\pi/2)$ satisfying $\theta,|\alpha|<\beta$. The double
periodicity condition amounts to $\frac{\sin \alpha}{\sin \beta}$ and
$\frac{\cos \alpha}{\cos \beta}$ being rational\footnote{take for 
instance $\beta \in (\pi/4,\pi/2)$ such that $\tan \beta$ is rational
and $\alpha = \pi/2-\beta$.}, so we write 
$\frac{\sin \alpha}{\sin \beta}=\frac{r}{s}$ and $\frac{\cos \alpha}{\cos 
\beta}
=\frac{p}{q}$. The lattice of periods $\Gamma$ is 
$\frac{q\pi}{\cos \beta}\Z \oplus i\frac{s\pi}{\sin\beta}\Z$, 
and the dual lattice is
\[
\Gamma^* = \frac{\cos\beta}{q\pi} \Z \oplus i\frac{\sin\beta}{s\pi}\Z
= \frac{\cos\alpha}{p\pi} \Z \oplus i\frac{\sin\alpha}{r\pi}\Z \; .
\]
Using the expression for the mean curvature vector\footnote{i.e. the half 
trace of the second fundamental form} is (in conformal coordinates) 
\[
H = \frac{e^{-2f}}{2} \left( \pa{\phi}{x} L_i 
\pa{X_{\theta,\beta}^\alpha}{x}
+ \pa{\phi}{y} L_i \pa{X_{\theta,\beta}^\alpha}{y}\right),
\]
where $\phi$ denotes here the Lagrangian angle, together with 
$\phi(z) = 2\pi\<\phi_0,z\>+\mathit{constant}\footnote{a careful 
computation 
shows that the constant is $\pi$.}$, we deduce that 
$\phi_0 = \frac{1}{\pi}\pa{\phi}{\bz} = \frac{e^{i\alpha}}{\pi}$ has 
lattice 
coordinates $(p,r)$ (in $\Gamma^*$). The periodicity condition above 
translates as the geometric property that the circle of radius 
$|\phi_0|=1/\pi$ possesses 8 lattice points (instead of the generic 4), 
namely: 
$\pm \frac{e^{i\alpha}}{\pi}, \pm \frac{e^{-i\alpha}}{\pi},
\pm \frac{e^{i\beta}}{\pi}, \pm \frac{e^{-i\beta}}{\pi}$. It may be 
that $\alpha = 0$, but the property still remains that there are 4 
extra points more than usual. These are exactly the points that come 
into play. Denoting $\gamma = \frac{e^{i\beta}}{2\pi} \in 
\frac{1}{2}\Gamma^*$,
\[
\Gamma^*_{\phi_0} = \{ \gamma,-\gamma,\bar{\gamma},-\bar{\gamma},
\frac{\bar{\phi}_0}{2}, -\frac{\bar{\phi}_0}{2} \}
\]
where the two last points are removed if $\phi_0$ is real (i.e. 
$\alpha = 0$). It also comes naturally that the limit case
$\alpha=\beta$ corresponds to the previous (and simpler) rectangular 
tori; if furthermore $\beta=\pi/4$ the lattice structure is exactly
that of the (square) torus.

The isometry described by I. Castro and
F. Urbano is $X_{\theta,\beta}^\alpha(z+it) = e^{t(\sin \alpha L_i
- \sin\beta R_i)} X(z)$. This property implies that the only
dual lattice elements in the Fourier exansion of $u=e^{\phi L_i/2}
\pa{X_{\theta,\beta}^\alpha}{z}$ are precisely 
$\gamma,-\gamma,\bar{\gamma}$
and $-\bar{\gamma}$; moreover opposite elements vanish simultaneously
and we have the conditions:
\[
\hat{a}_{-\gamma} = -e^{-i(\beta+\alpha)} \overline{\hat{a}_\gamma} \qquad
\hat{a}_{-\bar{\gamma}} = e^{i(\beta-\alpha)} \overline{\hat{a}_{\bar{\gamma}}}
\]
Using this the coefficients can be computed in terms of the functions
defined in~\cite{CU}.
\begin{figure}[h]
\begin{center}
\epsfig{file=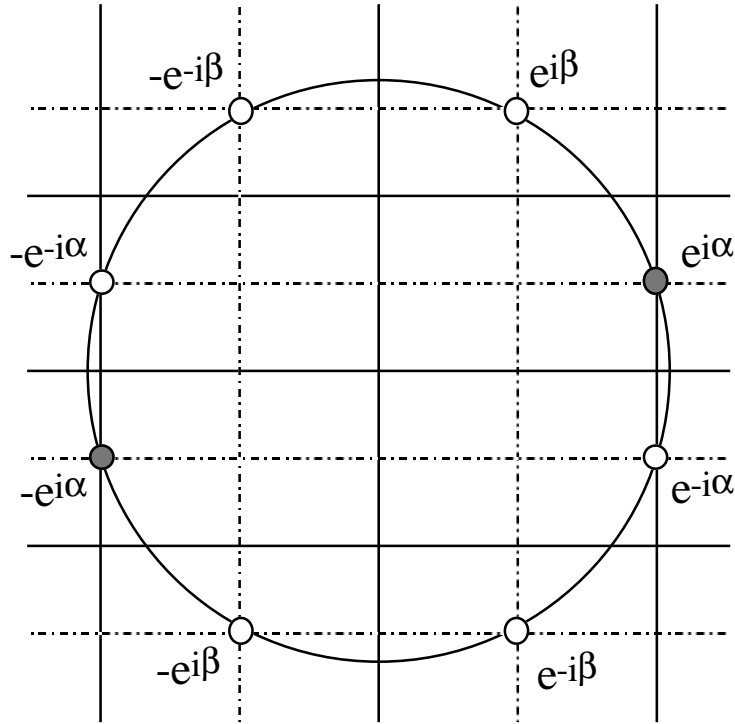, width=10cm}
\end{center}
\caption{A Castro and Urbano torus}
\end{figure}

%%%%%%%%%%%%%%%%%%%%%%%%%%%%%%%%%%%%%%%%%%%%%%%%%%%%%%%%%%%%%%%%%%%%%
%%%%%%%%%%%%%%%%%%%%%%%%%%%%%%%%%% LOOP GROUPS  %%%%%%%%%%%%%%%%%%%%%
%%%%%%%%%%%%%%%%%%%%%%%%%%%%%%%%%%%%%%%%%%%%%%%%%%%%%%%%%%%%%%%%%%%%%

\section{Introducing loop groups}

\subsection{Twisted loop groups}

We introduce loop groups, sets of maps $\lambda \mapsto G_{\lambda}$ from 
the circle $S^1 = \{\lambda \in \C/|\lambda |=1\}$
to some Lie groups (here various subgroups of ${\cal G}^{\C}$), with 
a multiplication law given as follows: the product of two elements
$\lambda \mapsto G_{\lambda}$ and $\lambda \mapsto G_{\lambda}'$ is just
$\lambda \mapsto G_{\lambda}G_{\lambda}'$. We denote

$$\Lambda {\cal G}:=\{[\lambda \mapsto G_{\lambda}]; S^1\rightarrow {\cal 
G}\}\hbox{ and }
\Lambda {\cal G}^{\C}:=\{[\lambda \mapsto G_{\lambda}]; 
S^1\rightarrow {\cal G}^{\C}\}.$$
We endow these groups with the $H^s$ topology for some $s>1/2$: if 
$G_{\lambda} = \sum _{k\in \Z}\hat{G}_k\lambda ^k$
is the Fourier expansion of $G_{\lambda}$, its $H^s$ norm is
$||G_{\lambda}||_s = \left(\sum _{k\in 
\Z}|\hat{G}_k|^2(1+k^2)^{s/2}\right) ^{1/2}$. Other topologies can be 
used (for instance
the ${\cal C}^{\infty}$ topology), for more details, see \cite{PS}.
We define the twisted loop groups 
\[
\Lambda \G_\tau = \{ [\lambda \mapsto G_{\lambda}] \in \Lambda {\cal G}/ 
G_{i\lambda} = \tau(G_\lambda) \}
\hbox{ and }\Lambda \G_\tau ^{\C}= \{ [\lambda \mapsto G_{\lambda}] 
\in \Lambda {\cal G}^{\C}/ G_{i\lambda} = \tau(G_\lambda) \},
\]
twisted meaning equivariant with respect to $\tau$. Also 
\[
\begin{array}{lll}
\Lambda^-_{\star} \GC_\tau &=& \{ [\lambda \mapsto G_{\lambda}] \in 
\Lambda \GC_\tau/ \textrm{$G_{\lambda}$ extends holomorphically}\\
&&\textrm{to the complement of the unit disk and $G_\infty = \1$} \}
\\
\Lambda^+ \GC_\tau &=& \{ [\lambda \mapsto G_{\lambda}] \in \Lambda 
\GC_\tau/
\textrm{$G_{\lambda}$ extends holomorphically to the unit disk} \}
\\
\Lambda^{+}_{\B}\GC_\tau &=& \{ [\lambda \mapsto G_{\lambda}] \in \Lambda 
\GC_\tau/
\textrm{$G_{\lambda}$ extends holomorphically to the unit disk and $G_0 
\in (\B,0)$} \}
\end{array}
\]
where $\B$ is some subgroup of ${\cal G}_0^{\C}$. In an analogous way 
define
the corresponding Lie algebras $\Lambda \g_\tau$,
$\Lambda \gC_\tau$, $\Lambda^-_\star \gC_\tau$, $\Lambda^+ \gC_\tau$
and $\Lambda_{\b}^+ \gC_\tau$ where \b\/ is the Lie algebra of $\B$.

$$
\begin{array}{lll}
\Lambda \gC_\tau &=& \{ [\lambda \mapsto \gamma _{\lambda}];S^1\rightarrow 
\g ^{\C}/ \gamma _{i\lambda}=\tau (\gamma _{\lambda}) \}
\\
\Lambda \g_\tau  &=& \{ [\lambda \mapsto \gamma _{\lambda}]\in \Lambda 
\gC_\tau /\gamma _{\lambda}\in \g ,\forall \lambda \in S^1 \}
\\
\Lambda^-_{\star} \gC_\tau &=& \{ [\lambda \mapsto \gamma _{\lambda}]\in 
\Lambda \gC_\tau /\gamma _{\lambda}\textrm{ extends holomorphically}\\
&&\textrm{to the complement of the unit disk and }\gamma_\infty = 0 \}
\\
\Lambda^+ \gC_\tau &=& \{ [\lambda \mapsto \gamma _{\lambda}]\in \Lambda_ 
\gC\tau /
\gamma _{\lambda}\textrm{ extends holomorphically to the unit disk} \}
\\
\Lambda^{+}_{\b} \gC_\tau &=& \{ [\lambda \mapsto \gamma _{\lambda}]\in 
\Lambda \gC_\tau /\gamma _{\lambda}
\textrm{ extends holomorphically to the unit disk and }\gamma _0\in (\b,0) 
\}
\end{array}
$$
An analysis of the relation $\gamma _{i\lambda}=\tau (\gamma _{\lambda})$, for
any $\gamma _{\lambda}\in \Lambda \gC_\tau$,
shows that, writing $\gamma _{\lambda}= \sum_{k\in\Z} \hat\gamma_k \lambda ^k$, this twisting condition is equivalent to
$\hat\gamma_k \in \gC_{k \bmod 4}$.

We remark in particular that $\Lambda \gC_\tau =  \Lambda^-_{\star} 
\gC_\tau \oplus \Lambda^+ \gC_\tau$, thus defining a
projection $[.]_{\Lambda^-_{\star} \gC_\tau}: \Lambda \gC_\tau \rightarrow 
\Lambda^-_{\star} \gC_\tau$.
Using this language, we can state the

\begin{corollary}       \label{extendedconnexion}
To each \g-valued 1-form $\alpha$ giving rise to a Hamiltonian stationary 
conformal 
Lagrangian immersion corresponds a $\Lambda \g_\tau$-valued 1-form 
$\alpha _\lambda$ (extended 1-form) satisfying relation (\ref{lambda}) and 
such that
        
\begin{equation}\label{k>-3}
        \left[ \alpha _{\lambda}\left( {\partial \over \partial z}\right) 
\right]_{\Lambda ^-_{\star}\gC _{\tau}} =
        \lambda ^{-2}\hat{\alpha}_{-2}\left( {\partial \over \partial 
z}\right) +
        \lambda ^{-1}\hat{\alpha}_{-1}\left( {\partial \over \partial 
z}\right) \hbox{ and }
        \left[ \alpha _{\lambda}\left( {\partial \over \partial 
\bar{z}}\right) \right]_{\Lambda ^-_{\star}\gC _{\tau}} =0,
\end{equation}
and
\begin{equation}\label{nonnul}
\hat{\alpha}_{-1}\left( {\partial \over \partial z}\right) \neq 0,
\end{equation}
and conversely. Moreover there exists a unique map $U_{\lambda}:\Omega 
\longrightarrow \Lambda {\cal G}_{\tau}$ such that
$dU_{\lambda} = U_{\lambda}\alpha _{\lambda}$ and $U_{\lambda}(z_0) = \1$. 
$U_{\lambda}$ is called an {\em extended lift}.
\end{corollary}
{\em Proof.} On one hand, Theorem 2 implies obviously that each $\g$-valued 
1-form $\alpha$ associated with a Hamiltonian
stationary conformal Lagrangian immersion can be deformed into such a 
$\alpha _{\lambda}$. On the other hand, any
$\Lambda \g _{\tau}$ valued 1-form $\alpha _{\lambda}$ satisfying 
(\ref{k>-3}) and (\ref{nonnul}) should satisfy

$$\alpha _{\lambda} = \lambda ^{-2}\hat{\alpha}_{-2} + \lambda 
^{-1}\hat{\alpha}_{-1}+ \hat{\alpha}_0
+ \lambda ^1\hat{\alpha}_1+ \lambda ^2\hat{\alpha}_2,$$
with $\hat{\alpha}_1({\partial \over \partial z}) = 
\hat{\alpha}_2({\partial \over \partial z}) = 0$ and
$\hat{\alpha}_1\left( {\partial \over \partial \bar{z}}\right) \neq 
0$, because of the reality condition
$\overline{\hat{\alpha}_k} = \hat{\alpha}_{-k}$ contained in the 
definition of $\Lambda \g _{\tau}$. If furthermore
$\alpha _{\lambda}$ satisfies (\ref{lambda}), then we conclude by using 
Theorem 2. The existence of $U_{\lambda}$ is
just a reformulation of Corollary 1.\bbox

\begin{remark} \em \lf
a) If $U$ is a LCLI and $U_{\lambda}$ is an extended LCLI such that $U_1 = 
U$, then the gauge action of
${\cal C}^{\infty}_{\star}(\Omega , {\cal G}_0)$ on $U$ extends in a 
natural way on $U_{\lambda}$. Precisely if
$K\in {\cal C}^{\infty}_{\star}(\Omega , {\cal G}_0)$
and if we denote $(KU)_{\lambda}$ the extended LCLI constructed from 
$U(K^{-1}, 0)$, then $(KU)_{\lambda} = U_{\lambda}(K^{-1}, 0)$.
To prove it, since we know that $(KU)_{\lambda}(z_0) = U_{\lambda}(K^{-1}, 
0)(z_0) = \1$, it suffices to check that both functions have
the same Maurer-Cartan form, namely

$$\lambda ^{-2}\alpha _2' + \lambda ^{-1}(K, 0)\alpha _{-1}' + (K, 
0)\alpha _0(K, 0)^{-1} - (dK.K^{-1},0) +
\lambda (K, 0)\alpha _1'' + \lambda ^2(K, 0)\alpha _2''.$$
b) The extended LCLI of the fundamental lift has the Maurer-Cartan form

$$\tilde{X}_{\lambda}^{-1}d\tilde{X}_{\lambda} = {\lambda ^{-2}\over 
2}{\partial \beta \over \partial z}(L_i,0)dz
+ \lambda ^{-1}e^f(0,\epsilon )dz + (M_X^{-1}dM_X, 0) + \lambda 
e^f(0,\bar{\epsilon})d\bar{z}
+ {\lambda ^2\over 2}{\partial \beta \over \partial 
\bar{z}}(L_i,0)d\bar{z},$$
which implies that $\tilde{X}_{\lambda} = (M_Xe^{\beta _{\lambda}L_i\over 
2}, X_{\lambda})$, where
$d\beta _{\lambda} = \lambda ^{-2}{\partial \beta \over \partial z}dz + 
\lambda ^2{\partial \beta \over \partial \bar{z}}d\bar{z}$.
Denoting $\gamma$ the harmonic conjugate function of $\beta$, i. e. such 
that ${1\over 2}(\beta +i\gamma )$ vanishes at $z_0$ and is
holomorphic, $ \beta _{\lambda} = {1\over 2}(\lambda ^{-2}+\lambda 
^2)\beta + {i\over 2}(\lambda ^{-2}-\lambda ^2)\gamma$.
The extended LCLI's of the spinors lifts have the Maurer-Cartan form

$$U_{\pm , \lambda}^{-1}dU_{\pm , \lambda} = {\lambda ^{-2}\over 
2}{\partial \beta \over \partial z}(L_i,0)dz
\pm \lambda ^{-1}(0,e^{-{\beta L_i\over 2}}{\partial X\over \partial z})dz
\pm \lambda (0,e^{-{\beta L_i\over 2}}{\partial X\over \partial 
\bar{z}})d\bar{z}
+ {\lambda ^2\over 2}{\partial \beta \over \partial 
\bar{z}}(L_i,0)d\bar{z},$$
and hence $U_{\pm , \lambda} = (\pm e^{\beta _{\lambda}L_i\over 2}, 
X_{\lambda})$.
\end{remark}

\subsection{Group decompositions}

The main tool for Weierstrass representations, as those proven in \cite{DPW}, 
are loop groups decompositions.
They are infinite dimensional analogs of Iwasawa decompositions such as 
$SU(n)^{\C} = SL(n,\C) = SU(n).\B$,
where $\B$ is a solvable (Borel) subgroup
of $SL(n,\C)$. For the convenience of the Reader, we first give here 
the proof of this splitting for the case $n=2$ (recall that
in our language, $SU(2)\simeq {\cal G}_0$).

\begin{proposition}
        Let $\B _0$ be the subgroup of matrices in ${\cal G}_0^{\C}$ 
leaving
        $\R^{*}_+\epsilon $ invariant, then ${\cal G}_0^{\C} = {\cal G}_0 
. \B _0$. More precisely the map
        $$\begin{array}{ccc}
        {\cal G}_0 \times \B _0 & \longrightarrow {\cal G}_0^{\C}\\
        (K,B) & \longmapsto KB
        \end{array}$$
        is a diffeomorphism.
\end{proposition}
\emph{Proof. } We use essentially Lemma~\ref{SU2action} and recall
that $\R^{\star}_+ \times {\cal G}_0$ acts freely and transitively on 
$\C \epsilon \oplus \C L_i \beps \setminus \{0\}$ which is the pointed $i$ 
eigenspace of $L_j$.
Since ${\cal G}_0$ commutes with $L_j$, so does ${\cal G}_0^{\C}$, hence 
$G\epsilon$ belongs to $\C \epsilon \oplus \C L_i \beps$ for any 
$G \in {\cal G}_0^{\C}$. By Lemma~\ref{SU2action}, there exist unique 
$K \in {\cal G}_0$ and $r\in \R ^{\star}_+$ such that $G\epsilon = 
rK\epsilon$.
Just set $B = K^{-1}G \in \B _0$. Notice that we might as 
well use $\bar{\B _0}$ to construct our Iwasawa decomposition. \bbox

Before stating the main results of this section, we shall establish a 
preliminary one. We set

$$\Lambda \vec{\cal G}_{\tau}^{\C} = \{ [\lambda \mapsto G_{\lambda}] 
\in \Lambda \vec{\cal G}^{\C}/
\tau (G_{\lambda}) = G_{i\lambda}, \forall \lambda \in S^1\}.$$
Notice that, since $\forall G\in \vec{\cal G}^{\C}$, $\tau ^2(G)=G$, 
any $G_{\lambda}\in \Lambda \vec{\cal G}_{\tau}^{\C}$
satisfies $G_{\lambda} = \tau ^2(G_{\lambda}) = G_{-\lambda}$. Also we 
denote

$$\Lambda {\cal G}_{2,\tau}^{\C} = \Lambda {\cal G}_2^{\C} \cap 
\Lambda {\cal G}_{\tau}^{\C}
= \{ [\lambda \mapsto K_{\lambda}] \in \Lambda {\cal G}_2^{\C}/ 
K_{\lambda}= \sum _{k\in \Z}\hat{K}_{2k}\lambda ^{2k},
\hat{K}_{4k}\in \C\1 ,\hat{K}_{4k+2}\in \C L_i\},$$
and

\begin{equation}\label{lambda4}
\Lambda {\cal G}_{0,\tau}^{\C} = \Lambda {\cal G}_0^{\C} \cap 
\Lambda {\cal G}_{\tau}^{\C}
= \{ [\lambda \mapsto F_{\lambda}] \in \Lambda {\cal G}_0^{\C}/ 
F_{\lambda}= \sum _{k\in \Z}\hat{F}_{4k}\lambda ^{4k} \}.
\end{equation}

\begin{lemma}\label{u1su2}
For any  $\lambda \mapsto G_{\lambda}\in \Lambda \vec{\cal 
G}_{\tau}^{\C}$, there exists
$(K_{\lambda},M_{\lambda})\in \Lambda {\cal G}_2^{\C}\times \Lambda 
{\cal G}_0^{\C}$,
unique up to sign, such that

$$G_{\lambda} = K_{\lambda}M_{\lambda}.$$
Moreover\\
(i) either $K_{\lambda}\in \Lambda {\cal G}_{2,\tau}^{\C}$ and 
$M_{\lambda}\in \Lambda {\cal G}_{0,\tau}^{\C}$,\\
(ii) or $K_{\lambda} = L_i\tilde{K}_{\lambda}$ and
$M_{\lambda} = \left( {1\over 2}(\lambda ^2+\lambda ^{-2})\1 +{1\over 
2i}(\lambda ^2-\lambda ^{-2})R_i\right) \tilde{M}_{\lambda}$,
with $(\tilde{K}_{\lambda},\tilde{M}_{\lambda})\in \Lambda {\cal 
G}_{2,\tau}^{\C}\times \Lambda {\cal G}_{0,\tau}^{\C}$.\\

In other words,
setting $\pi _{\lambda} :=  L_i\left( {1\over 2}(\lambda ^2+\lambda 
^{-2})\1 +{1\over 2i}(\lambda ^2-\lambda ^{-2})R_i\right)
\in \Lambda \vec{\G}_{\tau}$,

$$\Lambda \vec{\cal G}_{\tau}^{\C} = \Lambda {\cal 
G}_{2,\tau}^{\C}.\Lambda {\cal G}_{0,\tau}^{\C}\sqcup 
\pi _{\lambda}.\Lambda {\cal G}_{2,\tau}^{\C}.\Lambda {\cal 
G}_{0,\tau}^{\C}.$$
\end{lemma}
{\em Proof.} Let $\lambda \mapsto G_{\lambda}\in \Lambda \vec{\cal 
G}_{\tau}^{\C}$ and consider a lift
$g:\R \longrightarrow \vec{\cal G}^{\C}$ such that $g(\theta 
)=G_{e^{i\theta}}$, $\forall \theta \in \R $.
For any $\theta\in \R $, there exists $k(\theta )\in {\cal G}_2$ and 
$m(\theta )\in {\cal G}_0$ such that $g(\theta ) = k(\theta )m(\theta )$,
and $k(\theta )$ and $m(\theta )$ are unique up to sign. Moreover we can 
choose $k(\theta )$ and $m(\theta )$ to be continuous
functions of $\theta$. Since $G_{-\lambda} = G_{\lambda}$, we have 
$k(\theta +\pi )m(\theta +\pi ) = k(\theta )m(\theta )$,
and therefore $k(\theta +\pi ) = \pm k(\theta )$ and $m(\theta +\pi ) = 
\pm m(\theta )$. Hence
$k(\theta +2\pi ) = k(\theta )$ and $m(\theta +2\pi ) = m(\theta )$ and we 
may define
$(K_{\lambda},M_{\lambda})\in \Lambda {\cal G}_2^{\C}\times \Lambda 
{\cal G}_0^{\C}$
by $k(\theta )=K_{e^{i\theta}}$ and $m(\theta )=M_{e^{i\theta}}$. This 
proves the first assertion of the Lemma.

Notice that $\tau (K_{\lambda})\tau (M_{\lambda}) = \tau (G_{\lambda}) = 
G_{i\lambda} = K_{i\lambda}M_{i\lambda}$, which implies that
$(K_{i\lambda})^{-1}\tau (K_{\lambda}) = M_{i\lambda}(\tau 
(M_{\lambda}))^{-1} \in {\cal G}_2^{\C}\cap {\cal G}_0^{\C} = 
\{\pm \1 \}$.
Hence

\begin{equation}\label{subtwist}
\tau (K_{\lambda}) = s K_{i\lambda}\hbox{ and }\tau (M_{\lambda}) = s 
M_{i\lambda},\hbox{ with }s =\pm 1.
\end{equation}
Moreover because of the parity of $G_{\lambda}$, $K_{-\lambda}M_{-\lambda} 
= K_{\lambda}M_{\lambda}$, which leads to the alternatives\\

\noindent {\em a)} $K_{-\lambda}= K_{\lambda}$ and $M_{-\lambda}= 
M_{\lambda}$,\\
\noindent {\em b)} $K_{-\lambda}= -K_{\lambda}$ and $M_{-\lambda}= 
-M_{\lambda}$.\\

\noindent If {\em b)} occurs, $K_{\lambda}$ has the Fourier decomposition 
$K_{\lambda} = \sum _{k\in \Z}\hat{K}_{2k+1}\lambda ^{2k+1}$. Then
equation (\ref{subtwist}) implies that $\tau (\hat{K}_{2k+1})= s 
i(-1)^k\hat{K}_{2k+1}$, which is possible only if
all the $\hat{K}_{2k+1}$'s vanish, because the eigenvalues of the action 
of $\tau$ on $\C\1 +\C L_i$ are 1 and -1.
We exclude that since $K_{\lambda}\in {\cal G}_2^{\C}$. Hence only 
case {\em a)} may occur.

To conclude we inspect the consequence of (\ref{subtwist}). If $s =1$, 
case {\em (i)} of the Lemma occurs. If $s =-1$,
we define $\tilde{K}_{\lambda}$ in $\Lambda {\cal G}_2^{\C}$ and 
$\tilde{M}_{\lambda}$ in $\Lambda {\cal G}_0^{\C}$
by $K_{\lambda} = L_i\tilde{K}_{\lambda}$ and
$M_{\lambda} = \left( {1\over 2}(\lambda ^2+\lambda ^{-2})\1 +{1\over 
2i}(\lambda ^2-\lambda ^{-2})R_i\right) \tilde{M}_{\lambda}$.
Then we check that $\tau (\tilde{K}_{\lambda}) =\tilde{K}_{i\lambda}$ and 
$\tau (\tilde{M}_{\lambda}) =\tilde{M}_{i\lambda}$ which
shows that we are in case {\em (ii)}. \bbox

We recall results in \cite{PS}: let $\Gg$ be a compact Lie group
and $\Gg ^{\C}$ its complexification, and assume that the Iwasawa
decomposition
$\Gg ^{\C}=\Gg.{\cal B}_{\Gg}$ holds, for some solvable subgroup
${\cal B}_{\Gg}$ of $\Gg ^{\C}$. Define as before the loop groups
$\Lambda \Gg ^{\C}$, $\Lambda \Gg$, $\Lambda ^+\Gg ^{\C}$,
$\Lambda ^+_{{\cal B}_{\Gg}}\Gg ^{\C}$ and
$\Lambda ^-_{\star}\Gg ^{\C}$.

\begin{theorem}\label{ps}
[Pressley-Segal] \lf
a) The product mapping

$$\begin{array}{ccc}
\Lambda \Gg \times \Lambda ^+_{{\cal B}_{\Gg}}\Gg ^{\C} &
\longrightarrow & \Lambda \Gg ^{\C} \\
(\phi _{\lambda}, \beta _{\lambda}) & \longmapsto & \phi _{\lambda}. \beta
_{\lambda}
\end{array}$$
is a diffeomorphism.\\
b) There exists an open subset ${\cal C}_{\Gg}$ of $\Lambda \Gg
^{\C}$, called the {\em big cell}, such that the product mapping

$$\begin{array}{ccc}
\Lambda ^-_{\star}\Gg ^{\C}\times \Lambda ^+\Gg ^{\C} &
\longrightarrow & {\cal C}_{\Gg}\\
(\gamma ^-_{\lambda}, \gamma ^+_{\lambda}) & \longmapsto &
\gamma ^-_{\lambda}.\gamma ^+_{\lambda}
\end{array}$$
is a diffeomorphism.
\end{theorem}

We now use these results for proving the following decomposition theorems,
adapted to our
situation.

\begin{theorem}\label{thmA}
We have the decomposition
$\Lambda \GC _{\tau}= \Lambda \G _{\tau}.\Lambda ^+_{{\cal B}_0}\GC
_{\tau}$, i. e. the
map

$$\begin{array}{ccc}
\Lambda \G _{\tau}\times \Lambda ^+_{{\cal B}_0}\GC _{\tau} &
\longrightarrow &
\Lambda \GC _{\tau}\\
(F_{\lambda},B_{\lambda}) & \longmapsto & F_{\lambda}.B_{\lambda}
\end{array}$$
is a diffeomorphism.
\end{theorem}

\begin{theorem}\label{thmB}
There exists an open subset ${\cal C}$ of $\Lambda \G ^{\C}_{\tau}$,
called the {\em big cell}, such that
${\cal C} = \Lambda ^-_{\star}\G ^{\C}_{\tau}. \Lambda ^+\G
^{\C}_{\tau}$,
i. e. the product mapping

$$\begin{array}{ccc}
\Lambda ^-_{\star}\G ^{\C}_{\tau}\times \Lambda ^+\G 
^{\C}_{\tau} &
\longrightarrow & {\cal C}\\
(G^-_{\lambda}, G^+_{\lambda}) & \longmapsto & G^-_{\lambda}.G^+_{\lambda}
\end{array}$$
is a diffeomorphism.
\end{theorem}
{\em Proof of Theorem \ref{thmA}.} \lf
{\bf Step 1} We prove the
decomposition
$\Lambda \vec{\G}^{\C}_{\tau}=
\Lambda \vec{\G} _{\tau}.\Lambda ^+_{{\cal B}_0}\vec{\G}^{\C}_{\tau}$.

Let $G_{\lambda}\in \Lambda \vec{\G}^{\C}_{\tau}$. By Lemma 
\ref{u1su2},
$\exists (K_{\lambda}, M_{\lambda})\in \Lambda \G^{\C}_{2,\tau}
\times \Lambda \G^{\C}_{0,\tau}$ such that either
(i) $G_{\lambda} = K_{\lambda}. M_{\lambda}$, or
(ii) $G_{\lambda} = \pi _{\lambda}.K_{\lambda}. M_{\lambda}$.

We use Theorem \ref{ps} a) for $\Gg = SU(2)\simeq \G _0$. Let
$\tilde{M}_{\lambda} \in \Lambda \G^{\C}_0$ such that
$M_{\lambda} = \tilde{M}_{\lambda ^4}$. Then there exists a unique
$(\tilde{\phi}_{\lambda},\tilde{\beta}_{\lambda}) \in \Lambda \G_0 \times 
\Lambda ^+_{{\cal B}_0}\G^{\C}_0$
such that $\tilde{M}_{\lambda} = 
\tilde{\phi}_{\lambda}.\tilde{\beta}_{\lambda}$.
Setting $\phi _{\lambda} = \tilde{\phi}_{\lambda ^4}\in \Lambda 
\G_{0,\tau}$ (recall (\ref{lambda4})) and
$\beta _{\lambda} = \tilde{\beta}_{\lambda ^4}\in \Lambda ^+_{{\cal 
B}_0}\G_{0,\tau}^{\C}$,
we obtain $M_{\lambda} = \phi _{\lambda}\beta _{\lambda}$.

Similarly, we apply Theorem \ref{ps} a) for $\Gg = U(1)\simeq {\cal G}_2$:
since $K_{\lambda}\in \Lambda \G^{\C}_2$ there exists a unique
$(\psi _{\lambda}, \gamma _{\lambda})\in \Lambda \G_2 \times \Lambda 
^+_{{\cal B}_2}\G^{\C}_2$ such that
$K_{\lambda} = \psi _{\lambda}\gamma _{\lambda}$ (here we set
${\cal B}_2 = \{e^{itL_i}/t\in \R \}$.)
Thus $\tau (\psi _{\lambda}) \tau (\gamma _{\lambda}) = \tau (K_{\lambda}) 
= K_{i\lambda} = \psi _{i\lambda} \gamma _{i\lambda}$,
which implies that
$\tau (\psi _{\lambda})\psi _{i\lambda}^{-1} = \tau (\gamma 
_{\lambda})\gamma _{i\lambda}^{-1}
\in \Lambda \G _2\cap \Lambda ^+_{{\cal B}_2}\G_2^{\C} = \{\1 \}$.
(Here we used the fact that ${\cal B}_2$ is stable under the action of 
$\tau$
and therefore $\tau (\Lambda ^+_{{\cal B}_2}\G_2^{\C})\subset \Lambda 
^+_{{\cal B}_2}\G _2^{\C}$.)
Hence $\tau (\psi _{\lambda}) = \psi _{i\lambda}$ and
$\tau (\gamma _{\lambda}) = \gamma _{i\lambda}$, meaning that
$\psi _{\lambda}\in \Lambda \G _{2,\tau}$ and
$\gamma _{\lambda}\in \Lambda ^+_{{\cal B}_2}\G _{2,\tau}^{\C}$.
Lastly we remark that $\Lambda ^+_{{\cal B}_2}\G_{2,\tau}^{\C} =
\Lambda ^+_{\star}\G_{2,\tau}^{\C}$ and thus $\gamma _0 = \1$.

Hence we conclude that

$$G_{\lambda} = F_{\lambda} B_{\lambda},$$
where in case (i),

$$F_{\lambda} = \psi _{\lambda} \phi _{\lambda} \in  \Lambda \G 
_{2,\tau}.\Lambda \G_{0,\tau}\subset \Lambda \vec{\G}_{\tau},\hbox{ and }
B_{\lambda} = \gamma _{\lambda} \beta _{\lambda} \in \Lambda ^+_{\star}\G 
_{2,\tau}^{\C}.
\Lambda ^+_{{\cal B}_0}\G_{0,\tau}\subset \Lambda ^+_{{\cal 
B}_0}\vec{\G}_{\tau}.$$
And, in case (ii),

$$F_{\lambda} = \pi _{\lambda} \psi _{\lambda} \phi _{\lambda} \in \Lambda 
\vec{\G}_{\tau},\hbox{ and }
B_{\lambda} = \gamma _{\lambda} \beta _{\lambda} \in \Lambda ^+_{{\cal 
B}_0}\vec{\G}_{\tau}.$$
The diffeomorphism property of the decomposition is easy to check.\\

\noindent {\bf Step 2} We prove the decomposition
$\Lambda \G ^{\C}_{\tau} = \Lambda \G _{\tau}.\Lambda ^+_{{\cal 
B}_0}\G ^{\C}_{\tau}$.
Let $(G_{\lambda},T_{\lambda})\in \Lambda \G ^{\C}_{\tau}$. We want 
to prove
that there exist unique $(F_{\lambda},X_{\lambda})\in \Lambda \G _{\tau}$ 
and
$(B_{\lambda},b_{\lambda})\in \Lambda ^+_{{\cal B}_0}\G ^{\C}_{\tau}$,
such that

\begin{equation}\label{dcp1}
(F_{\lambda},X_{\lambda}) (B_{\lambda},b_{\lambda}) =
(F_{\lambda}B_{\lambda},F_{\lambda}b_{\lambda}+X_{\lambda}) =
(G_{\lambda},T_{\lambda}).
\end{equation}
Since $G_{\lambda}\in \Lambda \vec{\G}^{\C}_{\tau}$, the equation
$G_{\lambda} = F_{\lambda}B_{\lambda}$ has a unique solution
$(F_{\lambda},B_{\lambda})\in \Lambda \vec{\G}_{\tau}\cap \Lambda 
^+_{{\cal B}_0}\vec{\G}^{\C}_{\tau}$,
according to Step 1. The other equation,
$F_{\lambda}b_{\lambda}+X_{\lambda} = T_{\lambda}$, is equivalent to

\begin{equation}\label{dcpt1}
F^{-1}_{\lambda}T_{\lambda} = b_{\lambda} + F^{-1}_{\lambda}X_{\lambda}.
\end{equation}
Let us denote $\Lambda \C^4_{\tau} =  \{ [\lambda \mapsto 
V_{\lambda}];S^1\rightarrow \C^4/-L_jV_{\lambda}=V_{i\lambda} \}$,
$\Lambda \R ^4_{\tau} = \{ [\lambda \mapsto V_{\lambda}] 
;S^1\rightarrow \R[4]/ -L_jV_{\lambda}=V_{i\lambda} \}$
and $\Lambda ^+\C^4_{\tau} =  \{ [\lambda \mapsto V_{\lambda}]\in 
\Lambda \C^4_{\tau}/V_{\lambda}$ extends holomorphically to the unit 
disk $ \}$.
We have the following splitting

$$\Lambda \C^4_{\tau} = \Lambda \R ^4_{\tau} \oplus \Lambda 
^+\C^4_{\tau}.$$
We define $P:\Lambda \C^4_{\tau}\longrightarrow \Lambda 
\R^4_{\tau}$ to be
the projection on the first factor. Explicitely,

$$P\left (\sum _{n\in \Z}\hat{V}_{2n+1}\lambda ^{2n+1}\right) = \sum 
_{n\leq 0}\hat{V}_{2n-1}\lambda ^{2n-1} +
\sum _{n\geq 0}\overline{\hat{V}_{-2n-1}}\lambda ^{2n+1}.$$
Then the solution of (\ref{dcpt1}) is given by
$X_{\lambda} = F_{\lambda}P\left( F^{-1}_{\lambda}T_{\lambda}\right)$ and
$b_{\lambda} = F^{-1}_{\lambda}T_{\lambda} - P\left( 
F^{-1}_{\lambda}T_{\lambda}\right)$.\bbox
\noindent
\emph{Proof of Theorem \ref{thmB}}. We use Theorem (\ref{ps}) b) with
$\Gg = \vec{\G}\simeq U(2)$: let $\vec{\cal C} =
\Lambda ^-_{\star}\vec{\G}^{\C}.\Lambda ^+\vec{\G}^{\C}$ and
${\cal C} = \{(G_{\lambda},T_{\lambda})\in \Lambda \G ^{\C}_{\tau}/
G_{\lambda}\in \vec{\cal C}\}$.
The latter is clearly an open subset of $\Lambda \G ^{\C}_{\tau}$.
For any $(G_{\lambda},T_{\lambda})\in {\cal C}$, we look for 
$(G^-_{\lambda},T^-_{\lambda})\in \Lambda ^-_{\star}\GC _{\tau}$
and $(G^+_{\lambda},T^+_{\lambda})\in \Lambda ^+\GC _{\tau}$ such that
$(G^-_{\lambda}G^+_{\lambda}, G^-_{\lambda}T^+_{\lambda} + T^-_{\lambda}) 
= (G_{\lambda},T_{\lambda})$.
Let $(G^-_{\lambda}, G^+_{\lambda})$ be the unique element in $\Lambda 
^-_{\star}\vec{\G}^{\C}\times \Lambda ^+\vec{\G}^{\C}$
such that $G_{\lambda} = G^-_{\lambda} G^+_{\lambda}$. Since 
$G_{\lambda}\in \Lambda \vec{\G}^{\C}_{\tau}$,
$\tau (G^-_{\lambda}) \tau (G^+_{\lambda}) = \tau (G_{\lambda}) = 
G_{i\lambda} = G^-_{i\lambda}G^+_{i\lambda}$,
which implies
$$(G_{i\lambda})^{-1}\tau (G^-_{\lambda}) = G^+_{i\lambda}\tau 
(G^+_{\lambda})^{-1}\in 
\Lambda ^-_{\star}\vec{\G}^{\C}\cap \Lambda ^+\vec{\G}^{\C} 
=\{\1 \}.$$
Hence $(G^-_{\lambda}, G^+_{\lambda})\in \Lambda 
^-_{\star}\vec{\G}^{\C}_{\tau}\times \Lambda 
^+\vec{\G}^{\C}_{\tau}$.
To conclude, we look at the equation

\begin{equation}\label{dcpt2}
G^-_{\lambda}T^+_{\lambda} + T^-_{\lambda} = T_{\lambda}\Longleftrightarrow
T^+_{\lambda} + (G^-_{\lambda})^{-1}T^-_{\lambda} = 
(G^-_{\lambda})^{-1}T_{\lambda}.
\end{equation}
We let $\Lambda ^-\C^4_{\tau} =  \{ [\lambda \mapsto V_{\lambda}]\in 
\Lambda \C^4_{\tau}/V_{\lambda}$
extends holomorphically to the complement of the unit disk in $\C\cup 
\{\infty \} \}$ and use the linear splitting
$\Lambda \C^4_{\tau} = \Lambda ^-\C^4_{\tau}\oplus \Lambda 
^+\C^4_{\tau}$. Let
$Q^-: \Lambda \C^4_{\tau}\longrightarrow \Lambda ^-\C^4_{\tau}$ 
and
$Q^+: \Lambda \C^4_{\tau}\longrightarrow \Lambda ^+\C^4_{\tau}$
be the projection maps on each factor, namely
$Q^-\left( \sum _{n\in \Z}\hat{V}_{2n+1}\lambda ^{2n+1}\right) = \sum 
_{n\leq 0}\hat{V}_{2n-1}\lambda ^{2n-1}$
and
$Q^+\left( \sum _{n\in \Z}\hat{V}_{2n+1}\lambda ^{2n+1}\right) = \sum 
_{n\geq 0}\hat{V}_{2n+1}\lambda ^{2n+1}$.
Then the unique solution to (\ref{dcpt2}) is given by

$$T^-_{\lambda} = G^-_{\lambda}Q^-\left( 
(G^-_{\lambda})^{-1}T_{\lambda}\right) \hbox { and }
T^+_{\lambda} = Q^+\left( (G^-_{\lambda})^{-1}T_{\lambda}\right) .$$
Thus we obtained the right decomposition. \bbox

\section{Weierstrass representations}

\subsection{From Hamiltonian stationary surfaces to holomorphic potentials}

First we shall here sketch how to use ideas from \cite{DPW} in order to 
construct Weierstrass type data, starting from a Hamiltonian
stationary Lagrangian conformal immersion. Then we will revisit the 
obtained results and see how it simplifies in our situation.

Let $U = (F,X): \Omega \longrightarrow {\cal G}$ be a Hamiltonian 
stationary LCLI. Then it follows from Corollary~2, that
$U$ extends to a map $U_{\lambda} = (F_{\lambda},X_{\lambda}): \Omega 
\longrightarrow \Lambda {\cal G}_{\tau}$ satisfying
(\ref{k>-3}), (\ref{nonnul}) and $U_1 = U$.

\subsubsection{A family of holomorphic potentials}
There exists a holomorphic map $H_{\lambda} :\Omega \longrightarrow 
\Lambda {\cal G}_{\tau}^{\C}$ and a map
$B_{\lambda} :\Omega \longrightarrow \Lambda ^+_{{\cal B}_0}{\cal 
G}_{\tau}^{\C}$ such that

$$U_{\lambda}(z) = H_{\lambda}(z)B_{\lambda}(z),\ \forall \lambda \in S^1, 
\forall z\in \Omega .$$
The construction if $H_{\lambda}(z)$ and $B_{\lambda}(z)$ is done as 
follows: one looks for a map
$B_{\lambda} :\Omega \longrightarrow \Lambda ^+_{{\cal B}_0}{\cal 
G}_{\tau}^{\C}$ such that 
$H_{\lambda}(z) = U_{\lambda}(z)B_{\lambda}(z)^{-1}$ is holomorphic, i. e.

$$0 = {\partial (U_{\lambda}B_{\lambda}^{-1})\over \partial \bar{z}} 
= 
U_{\lambda}\left( \alpha _{\lambda}\left({\partial \over \partial 
\bar{z}}\right)  -
B_{\lambda}^{-1}{\partial B_{\lambda}\over \partial \bar{z}}\right)  
B_{\lambda}^{-1},$$
which is equivalent to

$${\partial B_{\lambda}\over \partial \bar{z}} = B_{\lambda}\left( 
\alpha _0\left({\partial \over \partial \bar{z}}\right) +
\lambda \alpha _1\left({\partial \over \partial \bar{z}}\right) +
\lambda ^2 \alpha _2\left({\partial \over \partial \bar{z}}\right) 
\right).$$
The existence of a solution $B_{\lambda}$ to this equation is first 
obtained locally (see \cite{DPW} or \cite{H2}), then
one can glue the local solutions into a global one \cite{DPW}. Then we write

$$H_{\lambda}^{-1} dH_{\lambda} = B_{\lambda}\left( \alpha _{\lambda} - 
B_{\lambda}^{-1}dB_{\lambda}\right) B_{\lambda}^{-1},$$
and using the fact that $B_{\lambda}$ takes its values in $\Lambda 
^+_{{\cal B}_0}{\cal G}_{\tau}^{\C}$ and that
$z\longmapsto H_{\lambda}(z)$ is holomorphic, we deduce that
\[
H_{\lambda}^{-1}dH_{\lambda} := \mu _{\lambda} = \sum _{n\geq 
-2}\hat{\mu}_n\lambda ^n,
\]
where each $\hat{\mu}_n$ is a closed (1,0)-form (i. e. holomorphic). As we 
shall see,  in  4.2, we can reconstruct
$U_{\lambda}$ from $\mu _{\lambda}$. J. Dorfmeister, F. Pedit and H.Y. Wu 
call the form $\mu _{\lambda}$ a {\em holomorphic potential}.
Notice that $\mu _{\lambda}$ is far from being uniquely defined, so we 
associate to $U_{\lambda}$ a whole
family of holomorphic potentials.

\subsubsection{A single meromorphic potential}

We can refine the above result as follows. First  
one can show that there exists a finite number of points
$a_1,...,a_k$ in $\Omega$ such that $U_{\lambda}(z)$ belongs to the big 
cell ${\cal C}$ (see Theorem \ref{thmB}), for all
$z\in \Omega \setminus \{a_1,...,a_k\}$. The proof of that is delicate and 
uses in particular the result of 5.1.1~.
Thus applying Theorem \ref{thmB}, we deduce that $\forall z\in \Omega 
\setminus \{ a_1,...,a_k\}$,
$\exists ! (U_{\lambda}^-(z),U_{\lambda}^+(z) )\in 
\Lambda ^-_{\star}{\cal G}_{\tau}^{\C}\times \Lambda ^+{\cal G}_{\tau}^{\C}$
such that

\begin{equation}\label{-+}
U_{\lambda}(z) = U_{\lambda}^-(z)U_{\lambda}^+(z),
\end{equation}
and then

\begin{equation}\label{mu}
\mu _{\lambda} := \left( U_{\lambda}^-\right) ^{-1}d U^-_{\lambda} = 
U_{\lambda}^+ \left( \alpha _{\lambda} - 
(U_{\lambda}^+)^{-1}dU_{\lambda}^+\right) \left( U_{\lambda}^+\right) 
^{-1}.
\end{equation}
We analyze equation (\ref{mu}): the right hand side tells us that 
$\hat{\mu}_n = 0$ for $n<-2$ and the
left hand side that $\hat{\mu}_n = 0$ for $n\geq 0$. Hence

\begin{equation}\label{mu1}
\mu _{\lambda} = \lambda ^{-2} \hat{\mu}_{-2} + \lambda ^{-1} 
\hat{\mu}_{-1}.
\end{equation}
Moreover, by writing the Fourier expansion of the right hand side of 
(\ref{mu}), one shows that
$\mu _{\lambda}\left({\partial \over \partial \bar{z}}\right) =0$. 
Hence $z\longmapsto U_{\lambda}^-(z)$ is holomorphic
on $\Omega \setminus \{a_1,...,a_k\}$. The analysis in \cite{DPW} shows 
furthermore that $z\longmapsto U_{\lambda}^-(z)$ extends
as a meromorphic map on $\Omega$: the potential $\mu _{\lambda}$ is a 
uniquely defined meromorphic potential.

\subsubsection{Explicit description}
We shall now revisit the previous facts. Since the 1-form $\mu _{\lambda}$ 
defined in (\ref{mu}) has his coefficients in
$\Lambda ^-_{\star}\g _{\tau}^{\C}$, we may write it as

\begin{equation}\label{mu2}
\mu _{\lambda} = \lambda ^{-2} (cL_i,0)dz + \lambda ^{-1} (0,a\epsilon 
+bL_i\bar{\epsilon})dz,
\end{equation}
where $a,b,c$ are meromorphic functions on $\Omega$.
Moreover, it follows from (\ref{mu}) that

$$(cL_i,0)dz = U^+_0 \alpha _2'\left( U^+_0\right) ^{-1} = 
U^+_0 ({1\over 2}{\partial \beta \over \partial z}L_i,0)\left( 
U^+_0\right) ^{-1}dz
=  ({1\over 2}{\partial \beta \over \partial z}L_i,0)dz,$$
where we used the fact that $U^+_0\in {\cal G}_0^{\C}$. Thus $c = 
{1\over 2}{\partial \beta \over \partial z}$.
Hence, letting $U_{\lambda}^- = (G_{\lambda}^-, T_{\lambda}^-)$ and using 
$dU_{\lambda}^- = U_{\lambda}^-\mu _{\lambda}$,
we obtain
\[
d(G_{\lambda}^-, T_{\lambda}^-) = \left( \lambda ^{-2} 
G_{\lambda}^-{1\over 2}{\partial \beta \over \partial z}L_idz,
\lambda ^{-1} G_{\lambda}^-(a\epsilon +bL_i\bar{\epsilon})dz\right) ,
\]
from which we deduce
\[
G_{\lambda}^-(z) = e^{{\lambda ^{-2}\over 4}(\beta (z)+i\gamma 
(z))L_i},
\]
where $\gamma :\Omega \longrightarrow \C$ 
is such that $\gamma (z_0) = 0$ and $d({1\over 2}(\beta +i\gamma )) = 
{\partial \beta \over \partial z}dz$
(\footnote{recall that ${1\over 2}(\beta +i\gamma )$ is the only 
holomorphic function vanishing at $z_0$ with a real part equal to ${\beta 
\over 2}$.}),
and also $dT_{\lambda}^- = \lambda ^{-1} e^{{\lambda ^{-2}\over 4}(\beta 
+i\gamma )L_i}(a\epsilon +bL_i\bar{\epsilon})dz$. Thus

$$(G_{\lambda}^-(z), T_{\lambda}^-(z)) = \left( e^{{\lambda ^{-2}\over 
4}(\beta (z)+i\gamma (z))L_i},
\lambda ^{-1}\int _{z_0}^z e^{{\lambda ^{-2}\over 4}(\beta (v)+i\gamma 
(v))L_i}(a(v)\epsilon +b(v)L_i\bar{\epsilon})dv \right) .$$
Now, letting  $U_{\lambda}^+ = (G_{\lambda}^+, T_{\lambda}^+)$,we can 
write (\ref{-+}) as

$$\left( e^{-{\lambda ^{-2}\over 4}(\beta (z)+i\gamma 
(z))L_i}F_{\lambda}(z), e^{-{\lambda ^{-2}\over 4}(\beta (z)+i\gamma 
(z))L_i}X_{\lambda}(z)\right) =$$
$$\left( G_{\lambda}^+(z), T_{\lambda}^+(z) +
e^{-{\lambda ^{-2}\over 4}(\beta (z)+i\gamma (z))L_i}\lambda ^{-1}\int 
_{z_0}^z e^{{\lambda ^{-2}\over 4}
(\beta (v)+i\gamma (v))L_i}(a(v)\epsilon +b(v)L_i\bar{\epsilon})dv 
\right) .$$
We conclude that 

$$\begin{array}{ccl}
(G_{\lambda}^-, T_{\lambda}^-) & = & \displaystyle \left(  e^{{\lambda 
^{-2}\over 4}(\beta +i\gamma )L_i} , 
e^{{\lambda ^{-2}\over 4}(\beta +i\gamma )L_i} Q^-\left( e^{-{\lambda 
^{-2}\over 4}(\beta +i\gamma )L_i}X_{\lambda}\right) \right) \\
(G_{\lambda}^+, T_{\lambda}^+) & = & \displaystyle \left( e^{-{\lambda 
^{-2}\over 4}(\beta +i\gamma )L_i}F_{\lambda} ,
Q^+\left( e^{-{\lambda ^{-2}\over 4}(\beta +i\gamma 
)L_i}X_{\lambda}\right) \right) .
\end{array}$$
After this analysis, we are led to the following.

\begin{theorem}
For any Hamiltonian stationary LCLI $U_{\lambda} = (F_{\lambda} 
,X_{\lambda} )$, there exist unique
$U^-_{\lambda} = (F^-_{\lambda} ,X^-_{\lambda} )\in \Lambda 
^-_{\star}{\cal G}_{\tau}^{\C}$ and
$U^+_{\lambda} = (F^+_{\lambda} ,X^+_{\lambda} )\in \Lambda ^+{\cal 
G}_{\tau}^{\C}$ such that
$U_{\lambda} = U^-_{\lambda}U^+_{\lambda}$, defined explicitely by

$$\begin{array}{cccl}
U^-_{\lambda} & = & \displaystyle \left(   e^{{\lambda ^{-2}\over 4}(\beta 
+i\gamma )L_i} ,
e^{{\lambda ^{-2}\over 4}(\beta +i\gamma )L_i} Q^-\left( e^{-{\lambda 
^{-2}\over 4}(\beta +i\gamma )L_i}X_{\lambda}\right) \right) \\
U^+_{\lambda} & = & \displaystyle \left(  e^{-{\lambda ^{-2}\over 4}(\beta 
+i\gamma )L_i}F_{\lambda} , 
Q^+\left( e^{-{\lambda ^{-2}\over 4}(\beta +i\gamma 
)L_i}X_{\lambda}\right) \right) ,
\end{array}$$
for $\gamma$ solution of $\gamma (z_0) = 0$ and $d ({1\over 2}(\beta 
+i\gamma ))= {\partial \beta \over \partial z}dz$.
Moreover,
\[
\mu _{\lambda} = \left( U^-_{\lambda}\right) ^{-1}dU^-_{\lambda} = 
\left( {\lambda ^{-2}\over 2}L_i{\partial \beta \over \partial z}dz, 
\lambda ^{-1}(a\epsilon +bL_i\bar{\epsilon})dz\right) 
\]
for some holomorphic functions $a,b$.
\end{theorem}
{\em Proof.} The uniqueness of the decomposition follows from Theorem 
\ref{thmB}. One checks easily that
$U_{\lambda} = U^-_{\lambda}U^+_{\lambda}$ and $U^-_{\lambda} \in \Lambda 
^-_{\star}{\cal G}^{\C}$.  For the verification of
$U^+_{\lambda} \in \Lambda ^+{\cal G}^{\C}$, we assume first that 
$U_{\lambda}$ corresponds to the fundamental lift:
then (see Remark 3) $F_{\lambda} = M_Xe^{\beta _{\lambda}L_i\over 2}$, 
which implies

$$G^+_{\lambda} = e^{-{\lambda ^{-2}\over 2}BL_i}F_{\lambda}
= M_Xe^{{\beta _{\lambda} - \lambda ^{-2}{1\over 2}(\beta +i\gamma )\over 
2}L_i}= M_Xe^{{\lambda ^2\over 4}(\beta -i\gamma )L_i}.$$
\noindent 
Therefore $G^+_{\lambda}$ belongs to $\Lambda ^+\vec{\cal G}^{\C}$. 
Thus obviously
$U^+_{\lambda} \in \Lambda ^+{\cal G}^{\C}$. If $U_{\lambda}$ 
corresponds to an arbitrary lift, then,
according to Remark 3, there exists $K\in {\cal C}^{\infty}(\Omega , {\cal 
G}_0)^{\star}$ such that
$F_{\lambda} = M_Xe^{\beta _{\lambda}L_i\over 2}K^{-1}$, and thus 
$G^+_{\lambda} =M_XK^{-1}e^{{\lambda ^2\over 4}(\beta -i\gamma )L_i}$
and we obtain the same conclusion.

Lastly repeating the argument of Theorem \ref{thmB}, we can deduce that 
$U^-_{\lambda} \in \Lambda ^-_{\star}{\cal G}_{\tau}^{\C}$
and $U^+_{\lambda} \in \Lambda ^+{\cal G}_{\tau}^{\C}$. The 
computation of $\mu _{\lambda}$ was done before. \bbox

The data $a,b,c={1\over 2}{\partial \beta \over \partial z}$ are called 
the {\em Weierstrass data of} $U_{\lambda}$.

\subsection{From a Weierstrass data to a Hamiltonian stationary conformal 
immersion}

We shall now see that the construction of the previous section has a 
converse. As above, we first sketch how to adapt the
strategy of \cite{DPW} and then we explore in more details what it means in our 
context.

Let $\mu _{\lambda} = \sum _{n\geq -2}\hat{\mu}_n\lambda ^n$ be a 
holomorphic potential; it is a 1-form on $\Omega$
with coefficients in $\Lambda \g _{\tau}^{\C}$ which is holomorphic, 
i. e. which satisfies
$\mu _{\lambda}\left( {\partial \over \partial \bar{z}}\right) = 0$ 
and $d\mu _{\lambda} = 0$. Then,
$\mu _{\lambda}\left( {\partial \over \partial \bar{z}}\right) = 0$ 
implies in particular that
$\mu _{\lambda}\wedge \mu _{\lambda} = 0$. Thus

$$d\mu _{\lambda} + \mu _{\lambda}\wedge \mu _{\lambda} = 0,$$
and there exists a unique map $H_{\lambda}\in \Lambda \G 
^{\C}_{\tau}$ such that $H_{\lambda}(z_0) = 0$ and
\begin{equation}\label{dH}
dH_{\lambda} = H_{\lambda}\mu _{\lambda}.
\end{equation}
For any $z\in \Omega$, we use Theorem \ref{thmA} with $H_{\lambda}(z)$: 
there exists a unique $(U_{\lambda}(z), V_{\lambda}(z))
\in \Lambda \G _{\tau}\times \Lambda ^+_{{\cal B}_0}\G ^{\C}_{\tau}$ 
such that $H_{\lambda}(z) = U_{\lambda}(z) V_{\lambda}(z)$.
A straightforward computation using (\ref{dH}) shows that

\begin{equation}\label{dF}
U_{\lambda}^{-1}dU_{\lambda} = V_{\lambda}\left( \mu _{\lambda} - 
V_{\lambda}^{-1}dV_{\lambda}\right) V_{\lambda}^{-1}.
\end{equation}
Let us denote $\alpha _{\lambda}:= U_{\lambda}^{-1}dU_{\lambda}$. 
Again the right hand side of (\ref{dF}) tells us that
$\alpha _{\lambda}$ should be of the form

$$\alpha _{\lambda} = \sum _{n\geq -2}\hat{\alpha}_n\lambda ^n,$$
but the left hand side says that $\overline{\hat{\alpha}_n} = 
\hat{\alpha}_{-n}$ and thus
\[
\alpha _{\lambda} = \lambda ^{-2}\hat{\alpha}_{-2} +  \lambda 
^{-1}\hat{\alpha}_{-1} +\hat{\alpha}_0
+\lambda \hat{\alpha}_1 + \lambda ^2 \hat{\alpha}_2.
\]
Moreover, using $V_{\lambda}^{-1} = \hat{V}^{-1}_0 +\lambda \hat{V}^{-1}_1 
+... = \hat{V}^{-1}_0 -\lambda \hat{V}^{-1}_0 \hat{V}_1\hat{V}^{-1}_0 
+...$,
it follows also from (\ref{dF}) that

$$\hat{\alpha}_{-2} = \hat{V}_0\hat{\mu}_{-2}\hat{V}^{-1}_0\hbox{ and }
\hat{\alpha}_{-1} = \hat{V}_0\hat{\mu}_{-1}\hat{V}^{-1}_0 + [\hat{V}_1, 
\hat{\mu}_{-2}]\hat{V}^{-1}_0$$
are (1,0)-forms. Hence, since $\alpha _{\lambda}$ satisfies condition 
(\ref{lambda}) automatically, Corollary 2 implies that -
provided that we can prove the condition $\hat{\alpha}_{-1}\neq 0$ - 
$U_{\lambda}$ is an extended lift of a Hamiltonian
stationary conformal immersion. Lastly, by the relation $U_{\lambda} = 
H_{\lambda}V_{\lambda}^{-1}$, we see that
$\mu_{\lambda}$ is a holomorphic potential for $U_{\lambda}$ in the sense 
of the above section.\\

Let us now look at the particular case where $\mu_{\lambda}$ has the form

$$\mu_{\lambda} = \left( {\lambda ^{-2} \over 2}{\partial \beta \over 
\partial z}L_i,\lambda ^{-1}(a\epsilon +bL_i\bar{\epsilon})
\right) dz,$$
for some holomorphic $\beta , a,b$. We integrate the equation 
$dH_{\lambda}= H_{\lambda}\mu _{\lambda}$. Denoting
$H_{\lambda} = (h_{\lambda}, \eta _{\lambda})$, it gives

$$(dh_{\lambda}, d\eta _{\lambda}) = \left( {\lambda ^{-2} \over 
2}{\partial \beta \over \partial z}h_{\lambda}L_idz,
\lambda ^{-1}h_{\lambda}(a\epsilon +bL_i\bar{\epsilon}) dz\right) .$$
It has the following solution

$$H_{\lambda}(z) = (h_{\lambda}(z), \eta _{\lambda}(z)) = \left( 
e^{{\lambda ^{-2} \over 4}(\beta (z)+i\gamma (z))L_i},
\int _{z_0}^z e^{{\lambda ^{-2} \over 4}(\beta (v)+i\gamma (v))L_i}\lambda 
^{-1}(a(v)\epsilon +b(v)L_i\bar{\epsilon}) dv\right) ,$$
where $\gamma$ is the harmonic conjugate function of $\beta$ vanishing at 
$z_0$. We now look for $U_{\lambda}=(F_{\lambda},X_{\lambda})
\in \Lambda \G _{\tau}$ and $V_{\lambda}=(B_{\lambda},b_{\lambda}) \in 
\Lambda ^+_{{\cal B}_0}\G _{\tau}^{\C}$ such that
$H_{\lambda} = U_{\lambda}V_{\lambda}$. We first use that $\beta 
_{\lambda} = {1\over 2}(\lambda ^{-2}+\lambda ^2)\beta
+ {i\over 2}(\lambda ^{-2}-\lambda ^2)\gamma$ (see the proof of Theorem 6) 
and thus

$$h_{\lambda} = e^{{\lambda ^{-2} \over 4}(\beta +i\gamma )L_i} = e^{{1 
\over 2}\beta _{\lambda}L_i}
e^{-{\lambda ^2 \over 4}(\beta -i\gamma )L_i},$$
meaning that we have $F_{\lambda} = e^{{1 \over 2}\beta _{\lambda}L_i}$ 
and $B_{\lambda} = e^{-{\lambda ^2 \over 4}(\beta -i\gamma )L_i}$.
Now we need to solve

$$\eta _{\lambda}(z) = \int _{z_0}^z e^{{\lambda ^{-2} \over 4}(\beta 
(v)+i\gamma (v))L_i}\lambda ^{-1}(a(v)\epsilon 
+b(v)L_i\bar{\epsilon})dv
= F_{\lambda}(z)b_{\lambda}(z) + X_{\lambda}(z),$$
or

$$e^{-{1 \over 2}\beta _{\lambda}(z)L_i}\int _{z_0}^z e^{{\lambda ^{-2} 
\over 4}(\beta (v)+i\gamma (v))L_i}
\lambda ^{-1}(a(v)\epsilon +b(v)L_i\bar{\epsilon})dv = b_{\lambda}(z) 
+e^{-{1 \over 2}\beta _{\lambda}(z)L_i}X_{\lambda}(z).$$

We deduce that

$$e^{-{1 \over 2}\beta _{\lambda}(z)L_i}X_{\lambda}(z) = P\left( e^{-{1 
\over 2}\beta _{\lambda}(z)L_i}
\int _{z_0}^z e^{{\lambda ^{-2} \over 4}(\beta (v)+i\gamma (v))L_i}\lambda 
^{-1}(a(v)\epsilon +b(v)L_i\bar{\epsilon})dv\right) .$$
Hence we proved 

\begin{theorem}
For any holomorphic datas $\beta , a,b$, the potential
$\mu_{\lambda} = \left( {\lambda ^{-2} \over 2}
{\partial \beta \over \partial z}L_i,
\lambda ^{-1}(a\epsilon +bL_i\bar{\epsilon}) \right) dz$
leads to construct the map $U_{\lambda}:\Omega \longrightarrow \Lambda
\G _{\tau}$ by
\[
U_{\lambda}(z) = (F_{\lambda}(z),X_{\lambda}(z)) =
\]\[
\left(
e^{{1 \over 2}\beta _{\lambda}(z)L_i}, e^{{1 \over 2}\beta _{\lambda}(z)L_i}
P\left( e^{-{1 \over 2}\beta _{\lambda}(z)L_i}
\int _{z_0}^z e^{{\lambda ^{-2} \over 4}(\beta (v)+i\gamma (v))L_i}\lambda ^{-1}(a(v)\epsilon +b(v)L_i\overline{\epsilon})dv\right) 
\right) ,$$
where $\beta _{\lambda} = {1\over 2}(\lambda ^{-2}+\lambda ^2)\beta
+ {i\over 2}(\lambda ^{-2}-\lambda ^2)\gamma$ and
$\gamma$ is the harmonic conjugate map to $\beta$ vanishing at $z_0$
(i.e. $\partial (\beta +i\gamma )/\partial \bar{z}=0$).
And $U_{\lambda}$ is an extended lift of a Hamiltonian stationary conformal
immersion if and only if $X_{\lambda}$ is an immersion.
\end{theorem}

%%%%%%%%%%%%%%%%%%%%%%%%%%%%%%%%%%%%%%%%%%%%%%%%%%%%%%%%%%%%%%%%%%%
%%%%%%%%   FINITE GAP SOLUTIONS    %%%%%%%%%%%%%%%%%%%%%%%%%%%%%%%%
%%%%%%%%%%%%%%%%%%%%%%%%%%%%%%%%%%%%%%%%%%%%%%%%%%%%%%%%%%%%%%%%%%%

\section{Tori and finite type solutions}

Going back to the torus, we will apply the concept of holomorphic 
potential defined in the previous section to the study of Hamiltonian 
stationary Lagrangian tori (in conformal coordinates). What makes the torus 
specific is that we can define -- intrinsically -- a notion of \emph{constant}
potential, i.e. $\mu = \eta dz$, where $dz$ is any globally defined holomorphic 
1-form and $\eta$ is a constant twisted loop of Lie-algebra elements. 
Indeed two globally defined holomorphic 1-forms on a torus \T\/ differ 
by a multiplicative constant. We may further restrict to those potentials 
having only a finite number of nonzero terms in their Fourier expansion
(known as \emph{polynomial loops}). While such conditions may seem 
(i) far-fetched and (ii) too restrictive, it turns out that 
\begin{itemize}
    \item  integrating potentials that are constant (in $z$) and 
    polynomial (in $\lambda$) is equivalent to integrating commuting flows, 
    which in our case leads to a much simpler integration process than 
    the Adler-Kostant-Symes (AKS) scheme; the corresponding HLC immersions are 
    called \emph{finite type solutions};

    \item  all immersed tori are finite type solutions.
\end{itemize}
Notice also that in the toric case, the considerations below prove 
the existence of potential (without resorting to the preceding section).
Such ideas originate in the theory of completely integrable systems,
however we will not explain here the link between commuting flows 
and finite type solutions, and refer the Reader to \cite{BFPP} 
for a good description of both sides of the AKS scheme.
Finally we will see how this new description relates to the one
given in section~\ref{linearproblem}.

\subsection{Construction of finite type solutions}

Throughout the section, $dz$ will denote some fixed global holomorphic 
1-form on a torus \T, or its universal cover \C. Then for any $d \in \N$ 
define
\[
\Lambda^d \g_\tau = \{ [\lambda \mapsto \xi _{\lambda}]\in \Lambda \g_\tau ; 
\xi_\lambda = \sum_{-d}^d \hat\xi_n \lambda^n \}
\]
the space of real polynomial loops of degree $d$.
\begin{proposition}
    Let $d \in 4\N+2$ and $\eta _{\lambda}\in \Lambda^d \g_\tau$ be a
polynomial loop. 
    Then the extended 1-form $\alpha_\lambda$ obtained through the AKS scheme 
    from the constant potential $\lambda^{d-2}\eta _{\lambda}dz$ on \C\/ (with starting 
    point $z_0$) is exactly the projection $\pi_{\Lambda \g_\tau}
    (\lambda^{d-2}\xi _{\lambda}dz)$ of the solution $\xi _{\lambda}$ to the following differential 
    equation:
    \begin{equation}	\label{lax1}
    \left\{ \begin{array}{l}
      d\xi _{\lambda}= [\xi _{\lambda}, \pi_{\Lambda \g_\tau}(\lambda^{d-2}
\xi _{\lambda}dz) ]
    \\
      \xi _{\lambda}(z_0) = \eta _{\lambda}
    \end{array} \right.
    \end{equation}
    where $\pi_{\Lambda \g_\tau}$ denotes the projection on 
    $\Lambda \g_\tau$ in the direct sum $\Lambda \gC_\tau = 
    \Lambda \g_\tau \oplus \Lambda^{+}_{\b} \gC_\tau$. Reciprocally,
    the solution exists and is complete.
\end{proposition}
\emph{Proof. }
First notice that $\lambda^{d-2}\eta _{\lambda}$ is a constant real polynomial loop
with lowest Fourier coefficient $\lambda^{-2}\hat{\eta}_{-d}$, thus also
a holomorphic potential, that we integrate on \C.
Let $M _{\lambda}\in \Lambda \GC_\tau$ be such that $M_\lambda(z_0) = \1$
and $\mu _{\lambda}= M_{\lambda}^{-1} dM_\lambda = 
\lambda^{d-2}\eta _{\lambda}dz$. Use the Iwasawa decomposition (Theorem \ref{thmA}) as 
in section 5
to write $M_\lambda = H_\lambda B_\lambda$ and by definition
$\alpha_\lambda = H_\lambda^{-1} dH_\lambda$. Set $\xi_\lambda (z)= 
H_\lambda^{-1} (z) \eta _{\lambda}H_\lambda$, which is well-defined on all \C.
By construction $\xi _{\lambda}$ is real (i.e. belongs to $\Lambda \G_\tau$)
since $\eta _{\lambda}$ and $H_\lambda (z)$ are. Using the fact that
$\eta _{\lambda}$ commutes with $M_\lambda$ we write
\[
\eta _{\lambda}= M_\lambda \eta _{\lambda}M_\lambda^{-1} 
= H_\lambda B_\lambda \eta _{\lambda}B_\lambda^{-1} H_\lambda^{-1}
\]
so $\xi_\lambda = B_\lambda \eta _{\lambda}B_\lambda^{-1}$, which proves that $\xi _{\lambda}$
has no Fourier coefficient with exponent lower than $-d$ (simply write
the Fourier expansions). Being real, $\xi _{\lambda}$ is $\Lambda^d \g_\tau$
valued. To prove that it solves the differential equation above, we write
$d\xi_\lambda = [\xi_\lambda , H_\lambda^{-1} dH_\lambda ]$; but
\begin{eqnarray*}
H_\lambda^{-1} dH_\lambda &=& 
\pi_{\Lambda \g_\tau} \left( H_\lambda^{-1} dH_\lambda \right)
=\pi_{\Lambda \g_\tau} \left( B_\lambda \mu _{\lambda}B_\lambda^{-1} 
- dB_\lambda B^{-1}_\lambda \right)
\\
&=& \pi_{\Lambda \g_\tau} \left( B_\lambda (\lambda^{d-2}\eta _{\lambda}dz) 
B_\lambda^{-1} \right)
=\pi_{\Lambda \g_\tau} \left( \lambda^{d-2} \xi_\lambda dz \right)
\end{eqnarray*}
{}\bbox

Any Hamiltonian stationary conformal Lagrangian immersion so obtained, 
either by integrating a constant polynomial loop as above, or by solving 
the differential 
equation~(\ref{lax1}), is called a \emph{finite type solution}.
Equation~(\ref{lax1}) can be written more explicitely, thus showing 
how to derive from $\xi_\lambda$ the extended 1-form $\alpha_\lambda$. 
Indeed writing $\xi _{\lambda}= \sum_{-d}^d \lambda^n \hat{\xi}_{n}$,
the projection $\pi_{\Lambda \g_\tau}(\lambda^{d-2} \xi _{\lambda}dz)$ is
\[
\pi_{\Lambda \g_\tau}(\lambda^{d-2} \xi _{\lambda}dz)
= \lambda^{-2} \hat{\xi}_{-d} dz + \lambda^{-1} \hat{\xi}_{-d+1} dz
+ \pi_{\g_0}(\hat{\xi}_{-d+2} dz)
+ \lambda \overline{\hat{\xi}_{-d+1}} d\bz
+ \lambda^{2} \overline{\hat{\xi}_{-d}} d\bz
\]
where $\pi_{\g_0}$ is the projection in the direct sum 
$\gC_0 = \g_0 \oplus \b$; define on $\gC_0$ the operator 
\[
r : \zeta \mapsto \frac{\pi_{\g_0}(\zeta)-i\pi_{\g_0}(i\zeta)}{2} 
\]
satisfying $\pi_{\g_0}(\zeta dz) = r(\zeta) dz + \overline{r(\zeta)}
d\bz$. Then we rewrite equation~(\ref{lax1}) as
\begin{equation}	\label{lax2}
\pa{\xi_\lambda}{z} = [ \xi_\lambda, \lambda^{-2} \hat{\xi}_{-d} 
+ \lambda^{-1} \hat{\xi}_{-d+1} + r(\hat\xi_{-d+2}) ]
\end{equation}
plus the initial condition; the conjugate equation is implied by the reality
of $\xi _{\lambda}$.
\lf

Since we aim at constructing solutions on a torus \T, we ought to 
notice that our construction, while valid on \C, does not necessarily 
give an immersion of the torus. To produce an actual torus we need 
to verify period conditions (a.k.a. monodromy conditions) obtained by 
integrating $\alpha_\lambda$. Also recall that the regularity of 
the immersion is equivalent to $\hat\xi_{-d+1}$ being non zero, 
otherwise the solution is only weakly conformal.

\subsection{A finiteness result}

\begin{theorem}
Let \T\/ be a 2-torus; then any Hamiltonian stationary conformal Lagrangian 
immersion in  \R[4] is of finite type.
\end{theorem}
This may seem surprising, especially if we think how restrictive the 
finite type condition is; however this result is almost classical 
in the theory of infinite dimensional integrable systems. As a consequence 
the space of solutions is a countable union of finite dimensional spaces.

\lf
\emph{Proof. }
We will adapt here an idea found in~\cite{BFPP}.
Let $X$ be a Hamiltonian stationary conformal Lagrangian immersion, 
$\alpha$ an associated Maurer-Cartan form (for some LCLI) and $\alpha_\lambda$ 
its extended 1-form. We first consider all quantities as being defined 
on the universal cover \C\/ of \T. We also choose a global
holomorphic 1-form $dz$ on \T\/ (and \C). A necessary and sufficient 
condition for $X$ to be of finite type is the existence of $\xi _{\lambda}: \C \to
\Lambda^d \g_\tau$ such that both equations hold
\begin{equation}	\label{lax3}
d\xi_\lambda = [ \xi_\lambda, \alpha_\lambda ]
\end{equation}
and
\begin{equation}	\label{adapted}
\alpha \left( \pa{}{z} \right) = \hat{\xi}_{-d} 
+ \hat{\xi}_{-d+1} + r(\hat\xi_{-d+2})
\end{equation}
Before we step into the proof, notice that finite type imposes a condition
obviously not satisfied in generality: consider the first Fourier
term (multiple of $\lambda^{-d}$) in equation~(\ref{lax3}), then
\[
d\hat\xi_{-d} = [ \hat\xi_{-d}, 
(r(\hat\xi_{-d+2}) - \hat\xi_{-d+2}) dz - \overline{\hat\xi_{-d+2}}d\bz ]
\in [ \g_2, \gC_0 ] = 0
\]
by the commutations properties of \g; so $\hat{\xi}_{-d}$ is constant.
Using condition~(\ref{adapted}), we see that $\hat{\xi}_{-d} = \frac{1}{2}
\pa{\beta}{z} (L_i,0)$, so that the Lagrangian angle is an affine function
of $x$ and $y$. On a torus this condition is trivially satisfied since
$\pa{\beta}{z}$ is holomorphic hence constant. This constant cannot vanish,
otherwise the immersed torus would be special-Lagrangian (or 
minimal); but there are no compact minimal tori.

The proof is divided in two steps: we first prove the existence
of a formal solution to~(\ref{lax3}) and (\ref{adapted}), then
extract a polynomial solution from these solutions,
using the property that each Fourier coefficient of the formal
solution satisfies an elliptic equation on the torus. By taking a
proper combination we infer the existence of $\xi _{\lambda}$.

\paragraph{Step 1: existence of adapted formal Killing fields. } 
A formal Killing field $\zeta _{\lambda}$ is a formal Laurent series in  
$\lambda$ verifying~(\ref{lax3}). Such a field is said \emph{adapted} 
if its first three terms are respectively equal to 
$\alpha_{-2}\left( \pa{}{z} \right)$, $\alpha_{-1}\left( \pa{}{z} \right)$ 
and $r\left(\alpha_0\left( \pa{}{z} \right) \right)$. 
Using the gauge action, we may suppose without loss of generality that 
$\alpha_0 = 0$, so that $\alpha \left( \pa{}{z} \right) = 
\alpha_{-2}\left( \pa{}{z} \right) + \alpha_{-1}\left( \pa{}{z} \right) = (aL_i,u)$ 
for some nonzero complex constant $a$, and $u(z) \in \C\epsilon \oplus 
\C L_i \beps$. Let us look for $\zeta _{\lambda}$ of the following form 
(typical of a gauge change used in~\cite{FT}):
$\zeta _{\lambda}= (\1,w_{\lambda})^{-1} (aL_i,u) (\1,w_{\lambda})$ where $w_{\lambda}$ has non negative
Fourier exponents (beware: $(\1,w_{\lambda})$ is a matrix in \GC, not in \gC).
The expression of $\zeta _{\lambda}$ simplifies here to give $\zeta _{\lambda}= (aL_i,u+aL_iw_{\lambda})$, 
and we solve the (1,0) part of equation~(\ref{lax3}):
\[
\left( 0 , \pa{u}{z} + aL_i \pa{w_{\lambda}}{z} \right) 
= [ ( aL_i,u+aL_iw_{\lambda}), (\lambda^{-2} aL_i, \lambda^{-1} u) ]
= (0, \lambda^{-1}aL_i u - \lambda^{-2} aL_i u + \lambda^{-2} a^2 w_{\lambda})
\]
so, using $a \neq 0$
\[
w_{\lambda} = \lambda^2 a^{-1} L_i \pa{w_{\lambda}}{z} + \lambda^2 a^{-2} \pa{u}{z} 
+ a^{-1} L_i u - \lambda a^{-1} L_i u
\]
Writing the hypothesis $w_{\lambda} = \sum_{n\geq 0} \lambda^n \hat{w}_n$, we obtain
$w_{\lambda}$ by simple recurrence (hence the formal series):
\[
\left\{ \begin{array}{l}
\hat{w}_0 = a^{-1}L_i u	\\
\hat{w}_1 = -a^{-1} L_i u  \\
\\
\ds \hat{w}_2 = a^{-1}L_i \pa{\hat{w}_0}{z} + a^{-2} \pa{u}{z} = 0	\\
\\
\ds \hat{w}_n = a^{-1}L_i \pa{\hat{w}_{n-2}}{z} \quad \textrm{for $n>2$}
\end{array} \right.
\]
So $w_{\lambda} = a^{-1}L_i u - \sum_{n \geq 0} \lambda^{2n+1} (a^{-1}L_i)^{n+1} 
\pa[n]{u}{z}$. Finally $\zeta _{\lambda}= (aL_i, \sum_{n \geq 0} 
\lambda^{2n+1} a^{-n}L_i^n \pa[n]{u}{z})$. We now check that the (0,1)
equation holds. Using the same idea as in section~\ref{linearproblem},
$d\alpha'_{-1} + [\alpha''_1 \wedge \alpha'_2]+[\alpha'_{-1} \wedge \alpha''_0]
=0 $, which in our notations yields $\pa{u}{\bz} = aL_i \bar{u}$. Then
\[
\pa{\zeta_\lambda}{\bz} = \left( 0, 
\sum_{n \geq 0} \lambda^{2n+1} a^{-n}L_i^n 
\frac{\partial^{n+1}u}{\partial \bz \partial z^n} \right)
=\left( 0, aL_i\lambda \bar{u} - \lambda^2 \bar{a} L_i \sum_{n \geq 0} 
\lambda^{2n+1} a^{-n}L_i^n \pa[n]{u}{z}\right)
= [\zeta_\lambda, (\lambda^2 \bar{a} L_i, \lambda \bar{u})] 
\]
Finally we verify that $\zeta _{\lambda}$ is adapted: $\hat{\zeta}_0 = (aL_i,0) 
= \alpha'_{-2}\left( \pa{}{z} \right)$,
$\hat{\zeta}_1 = (0,u) = \alpha'_{-1}\left( \pa{}{z} \right)$ and
$\hat{\zeta}_2 = 0$. 
It should be noted that for any $n \in \Z$, $\lambda^n\zeta_\lambda$ is still 
an adapted formal Killing field.

\paragraph{Step 2: elliptic equation and polynomial Killing fields. }

The possibility of reducing formal Killing fields to polynomial ones
relies on the following property: all coefficients of the formal series 
$\zeta _{\lambda}$ satisfy the elliptic equation
\begin{equation} 	\label{ellipticeqn}
\left( \pa{}{z} \pa{}{\bz} +|a|^2 \right) \hat{\zeta}_n = 0
\end{equation}
Recall that $u$ being doubly periodic (since originally defined on \T), 
so are all these coefficients.
The space of solutions to an elliptic equation on a compact Riemann surface
is finite dimensional. Consider now the free sequence $(\zeta^m_{\lambda})_{m \geq 0}$ 
of truncated formal Killing fields: $\zeta^m_\lambda = 
(\lambda^{-4m-2}\zeta_\lambda)_-$, truncated meaning that we keep only
negative powers of $\lambda$. The image of the sequence by the operator
$d+\ad_{\alpha_\lambda}$ is finite dimensional, indeed for any 
$\phi _{\lambda}= \lambda^{-4m-2}\zeta _{\lambda}= \sum_{n \geq -4m-2} \lambda^n \hat\phi_n$
\begin{eqnarray*}
(d+\ad_{\alpha_\lambda})(\phi _{\lambda})_- &=& d(\phi _{\lambda})_- + [\alpha_\lambda,(\phi _{\lambda})_-] 
= (d\phi _{\lambda})_- - [(\phi _{\lambda})_-,\alpha_\lambda]
= [\phi _{\lambda},\alpha_\lambda]_- - [(\phi _{\lambda})_-,\alpha_\lambda]
\\
&=& [\lambda^{-2} \hat\phi_{-2} + \lambda^{-1} \hat\phi_{-1}
+ \hat\phi_0 + \lambda \hat\phi_1 + \lambda^2 \hat\phi_2, \alpha_\lambda]_-
 - [\lambda^{-2} \hat\phi_{-2} + \lambda^{-1} \hat\phi_{-1}
+ \hat\phi_0, \alpha_\lambda]
\end{eqnarray*}
All other terms either vanish in the truncature or compensate between 
the two brackets. Since each coefficient $\hat\phi_k$ belongs to a
finite dimensional space, we have our claim. So a finite combination of the 
$\zeta^m_{\lambda}$, call it $\phi$, lies in the kernel of $d+\ad_{\alpha_\lambda}$,
and is automatically adapted. Then $\xi _{\lambda}= \phi _{\lambda}+ \bar{\phi_{\lambda}}$ is adapted,
satisfies equation~(\ref{lax3}) and is real; that is $\xi _{\lambda}$
belongs to some $\Lambda^d \g_\tau$. \bbox

\begin{remark} \em 	\label{remevencoef}
The solutions of~(\ref{lax3}) we have constructed have the following
property: $\xi_\lambda$ has no term of Fourier exponent equal to
$0 \pmod 4$, except the first and last ones; $\hat{\xi}_n = 0$ but for
$\hat{\xi}_{-d} = (aL_i,0)$ and $\hat{\xi}_d=\overline{\hat{\xi}_{-d}}$.
\end{remark}

\begin{remark} \em
So far we have used a fixed complex coordinate $z$ on \T\/ (or \C),
but we might want to switch to another coordinate say $w = \mu z$. A quick 
look at~(\ref{lax1}) shows that a solution in the $z$ variable
is usually not valid in the $w$ variable (that can also be seen 
on~(\ref{lax2})). However the immersion stays of finite type whatever
the coordinate may be; so an Hamiltonian stationary immersion can be of 
finite type $d$ for some variable $z$ and $d' \neq d$ for another variable $w$.
The type itself is not a well-defined invariant, and an example of this will be 
given in the next section.
\end{remark}

\subsection{Finite type and lattice properties}

We will now use the information given by section~\ref{linearproblem}
together with the finite type point of view. Let $\Gamma$ be a lattice, 
with dual lattice $\Gamma^*$, and set $\T = \C/\Gamma$. For any Hamiltonian 
stationary conformal Lagrangian immersion $X$ there exists $p \in \N$ 
and $\xi_\lambda$ in $\Lambda^{4p+2} \g_\tau$ solution of~(\ref{lax3})
projecting to the extended 1-form $\alpha_\lambda$ associated with $X$ 
(more precisely with one of the spinor lifts of $X$). As mentioned in 
remark~\ref{remevencoef} above, we may assume that $\xi_\lambda$ 
has no term of Fourier exponent equal to $0 \bmod 4$. So we can write
\[
\xi_\lambda = \lambda^{-4p-2}(\frac{\pi\bar{\beta}_0}{2}L_i,0)
+\sum_{q=-p}^{p} \left( \lambda^{4q-1} (0,u_q) + \lambda^{4q+1} (0,v_q) 
+ \lambda^{4q+2} (\pi \frac{\pi\bar{\beta}_0 c_q}{2} L_i,0) \right)
\]
with $\beta_0 \in \Gamma^{*} - {0}$, $u_q(z) \in \gC_{-1}$ and
$v_q(z) \in \gC_{1}$. 
Equation~(\ref{lax2}) can be rewritten using the Fourier expansion, 
\begin{equation}	\label{recfourier}
\pa{\hat{\xi}_n}{z} = [\hat{\xi}_{n+2}, \hat{\xi}_{-4p-2}]
+ [\hat{\xi}_{n+1}, \hat{\xi}_{-4p-1}] .
\end{equation}
Commutation properties $[\gC_0,\gC_2]=[\gC_2,\gC_2]=[\gC_1,\gC_{-1}]=0$
show that terms with even Fourier exponent -- hence the $c_{q}$'s --
are constant (and given by 
the initial condition in~(\ref{lax1})). We can then use~(\ref{recfourier}) 
or its conjugate to derive a recurrence relation:
\[
v_q = \frac{2}{\pi \bar{\beta}_0} L_i \pa{u_q}{z}
\]
\[
u_{q+1} = c_q u_{-p} + \frac{2}{\pi \bar{\beta}_0} L_i \pa{v_q}{z}
= c_q u_{-p} - \left( \frac{2}{\pi \bar{\beta}_0}\right)^2  \pa[2]{u_q}{z}\; .
\]
Conjugating equation~(\ref{recfourier}) (or recalling~(\ref{ellipticeqn})) 
one easily derives the second order equation 
$\Delta u_q + \pi^2 |\beta_0|^2 u_q = 0$, so all terms have Fourier frequencies 
$\gamma \in \frac{1}{2}\Gamma^*$ such that $|\gamma| = \frac{1}{2} |\beta_0|$. 
Thus we write $u_q = \sum_{\gamma} (\hat{a}_{\gamma,q} \epsilon
+ \hat{b}_{\gamma,q} L_i \bar{\epsilon}) e^{2i\pi\< \gamma,z\>}$
and the recurrence relation yields
\[
\hat{a}_{\gamma,q+1} = c_q \hat{a}_{\gamma,-p} + 
\left( \frac{2\bar{\gamma}}{\bar{\beta}_0} \right)^2 \hat{a}_{\gamma,q}
\]
(same equation for $\hat{b}_{\gamma,q}$). Taking the first and last
$\gC_{-1}$ terms:
\[
\hat{a}_{\gamma,p} = \hat{a}_{\gamma,-p} \left( \sum_{q=0}^{2p}
c_{p-q-1} \left( \frac{2\bar{\gamma}}{\bar{\beta}_0} \right)^{2q} 
\right)
\]
with the convention that $c_{-p-1} = 1$. Now we may compare
both ends of the chain using the reality of $\xi_\lambda$: 
$v_p = \bar{u}_{-p}$, while
\[
v_p = \frac{2}{\pi \bar{\beta}_0} L_i \pa{u_p}{z}
\]
which yields (recalling from section~\ref{linearproblem} the 
expression for $u_{-p}$): $\gamma\hat{a}_{\gamma,-p} 
= -\bar{\gamma} \hat{a}_{\gamma,p}$. So for each frequency $\gamma$
such that $\hat{a}_{\gamma,-p} \neq 0$
\[
\sum_{q=0}^{2p} c_{p-q-1} 
\left( \frac{2\bar{\gamma}}{\bar{\beta}_0} \right)^{2q} 
= -\, \frac{\gamma}{\bar{\gamma}} \; .
\]
We may rewrite this condition as a polynomial equation of degree exactly 
$d=4p+2$ in $\gamma$:
\begin{equation}	\label{polynome}
\gamma^{4p+2} + \frac{\bar{\beta}_0}{2} 
\left( \frac{\beta_0}{2} \right)^{4p+1} \sum_{q=0}^{2p} 
c_{q-p-1} \left( \frac{2\gamma}{\beta_0}\right)^{2q}=0  \; .
\end{equation}
We conclude with the following result:
\begin{theorem}
A finite type solution with type $d \in 4\N+2$ and Lagrangian angle 
$\beta = 2\pi\<\beta_0,z\> + \mathit{constant}$ has all its Fourier 
frequencies in $\Gamma^*_{\beta_0}$ and satisfying a polynomial equation 
of degree $d$, depending only on $\beta_0$ and the initial value 
in~(\ref{lax3}) (i.e. the constant potential). As a consequence 
$\mathrm{Card}(\Gamma^*_{\beta_0}) \leq d$.
\end{theorem}

\paragraph{Genus zero solutions. }
We conclude with the study of the simplest case of type $d=2$,
also called genus zero solutions. Then condition~(\ref{polynome}) implies 
that $\gamma = \pm i |\frac{1}{2}\beta_0|$. There are only two
possibilities, and if the lattice is rectangular, we find -- after a change 
of variable -- the rectangular generalizations of the standard torus.
It should be noted that the condition on $\gamma$ implies that
\emph{-- in the original complex coordinate --} the term
$e^{-\beta L_i/2} \pa{X}{z}$ depends only on $y$. This is clearly
variable-dependant, and shows that the type $d$ may well change
if one changes the variable.

%%%%%%%%%%%%%%%%%%%%%%%%%%%%
%%%%%%%%%%  REFERENCES%%%%%%
%%%%%%%%%%%%%%%%%%%%%%%%%%%%

\newpage %%%%%%%%%%%%%%%%%

\lf\lf
{\small 
Fr\'ed\'eric H\'elein \lf
CMLA, ENS de Cachan, 61 avenue du Pr\'esident 
Wilson, 94235 Cachan Cedex, France \lf
\texttt{helein@cmla.ens-cachan.fr} \lf
\lf
Pascal Romon \lf
CMLA, ENS de Cachan, 61 avenue du Pr\'esident Wilson, 94235 Cachan 
Cedex, France, and Universit\'e de Marne-la-Vall\'ee, France \lf
\texttt{romon@cmla.ens-cachan.fr}
}

\end{document}